\documentclass{amsart}

\usepackage{latexsym,enumerate}
\usepackage{amsmath,amsthm,amsopn,amstext,amscd,amsfonts,amssymb}
\usepackage[ansinew]{inputenc}
\usepackage{verbatim}
\usepackage{epsfig} 
\usepackage{graphicx,psfrag} 
\usepackage{subfigure}

\newcommand{\CC}{\mathbb{C}}

\newcommand{\RR}{\mathbb{R}}
\newcommand{\ZZ}{\mathbb{Z}}

\newcommand{\DD}{\mathbb{D}}

\newtheorem{teorema}{Theorem}[section]
\newtheorem{lema}[teorema]{Lemma}
\newtheorem{prop}[teorema]{Proposition}
\newtheorem{corolario}[teorema]{Corollary}
\newtheorem{definicion}[teorema]{Definition}

\DeclareMathOperator{\cotanh}{cotanh}

\DeclareMathOperator{\Arcsinh}{Arcsinh}
\DeclareMathOperator{\Arccosh}{Arccosh}

\newcommand{\spb}[1]{\smallskip}
\newcommand{\mpb}[1]{\medskip}
\newcommand{\bpb}[1]{\bigskip}

\newcommand{\p}{\partial}

\renewcommand{\d}{\delta}
\newcommand{\D}{\Delta}
\newcommand{\g}{\gamma}

\renewcommand{\l}{\lambda}

\renewcommand{\O}{\Omega}
\newcommand{\s}{\sigma}

\begin{document}

\title{A REAL VARIABLE CHARACTERIZATION OF GROMOV HYPERBOLICITY OF FLUTE SURFACES}

\author{
Ana Portilla$^{(1)}$, Jos\'e M. Rodr{\'\i}guez
$^{(1)}$ \and Eva Tour{\'\i}s$^{(1)}$}
\date{April 11, 2008.\\
$(1)\,\,\,$ Research partially supported by three grants from
M.E.C. (MTM 2006-11976, MTM 2006-13000-C03-02 and MTM 2007-30904-E), Spain.}

\maketitle{}



Departamento de Matem\'aticas

Escuela Polit\'ecnica Superior

Universidad Carlos III de Madrid

Avenida de la Universidad, 30

28911 Legan\'es (Madrid), SPAIN



emails: {\tt apferrei@math.uc3m.es, jomaro@math.uc3m.es, etouris@math.uc3m.es}

Corresponding author: Jos\'e M. Rodr{\'\i}guez

\begin{abstract}
In this paper we give a characterization of the Gromov
hyperbolicity of trains (a large class of Denjoy domains which
contains the flute surfaces) in terms of the behavior of a real
function. This function describes somehow the distances between
some remarkable geodesics in the train. This theorem has several
consequences; in particular, it allows to deduce a result about
stability of hyperbolicity, even though the original surface and
the modified one are not quasi-isometric.
\end{abstract}

\


\spb


\bigskip

\section{Introduction.}

\mpb

The theory of Gromov hyperbolic spaces is a useful tool in order to understand the connections between graphs and
Potential Theory (see e.g. \cite{ARY}, \cite{CFPR}, \cite{FR2},
\cite{HS}, \cite{K1}, \cite{K2}, \cite{K3}, \cite{R1}, \cite{R2},
\cite{So}). Besides, the
concept of Gromov hyperbolicity grasps the essence of negatively
curved spaces, and has been successfully used in the theory of
groups (see e.g. \cite{GH}, \cite{G1}, \cite{G2} and the references
therein).

A geodesic metric space is called hyperbolic (in the Gromov sense)
if there exists an upper bound of the  distance of every point in a
side of any geodesic triangle to the union of the two other sides
(see Definition \ref{def:Rips}). The latter condition is known as
Rips condition.

But, it is not easy to determine whether a given space is Gromov
hyperbolic or not.
Recently, there has been some research aimed to show that metrics
used in geometric function theory are Gromov hyperbolic.
Some specific examples are showing that the Klein-Hilbert metric
(\cite{Be}, \cite{KN}) is Gromov hyperbolic
(under particular conditions on the domain of definition),
that the Gehring-Osgood metric
(\cite{Ha}) is Gromov hyperbolic,
and that the Vuorinen metric
(\cite{Ha}) is not Gromov hyperbolic
(except for a particular case).
Recently, some interesting results
by Balogh and Buckley \cite{BB}
about the hyperbolicity
of Euclidean bounded domains
with their quasihyperbolic metric
have made significant progress in this direction
(see also \cite{BHK}, \cite{V} and the references therein).
Another interesting instance is that of a Riemann
surface endowed with the Poincar\'e metric. With such metric
structure a Riemann surface is always negatively curved, but not every
Riemann surface is Gromov hyperbolic, since topological obstacles
may impede it: for instance, the two-dimensional jungle-gym (a
$\ZZ^2$-covering of a torus with genus two) is not hyperbolic.

We are interested in studying when Riemann surfaces equipped with
their Poincar\'e metric are Gromov hyperbolic
(see e.g. \cite{RT1}, \cite{RT2}, \cite{RT3},
\cite{PRT1}, \cite{PRT2}, \cite{PRT3}, \cite{APRT}, \cite{PT}).
To be more precise, in the current paper
our main aim is to study the hyperbolicity of Denjoy domains, that
is to say, plane domains $\O$ with $\partial \O \subset \RR$. This
kind of surfaces are becoming more and more important in Geometric
Theory of Functions, since, on the one hand, they are a very general
type of Riemann surfaces, and, on the other hand, they are more
manageable due to its symmetry.
For instance, Garnett and Jones have proved the Corona Theorem for Denjoy
domains (\cite{GJ}), and in \cite{APR} the authors
have got the characterization of Denjoy domains
which satisfy a linear isoperimetric inequiality.

Denjoy domains are such a wide class of Riemann surfaces that
characterization criteria are not straightforward to apply. That is
the main reason that led us to focus on a particular type of Denjoy
domain, which we have called \emph{train}. A train can be defined as the complement of a sequence of ordered closed intervals (see Definition
\ref{def:Train}).
Trains do include a especially important case of surfaces which are
the flute surfaces (see, e.g. \cite{B1}, \cite{B2}).
These ones are the simplest examples of
infinite ends, and besides, in a flute surface it is possible to give a fairly
precise description of the ending geometry (see, e.g. \cite{H}).
In \cite{APRT} there are some partial results on hyperbolicity of trains.

This paper is a natural continuation of \cite{APRT}. Although some of the theorems in the current work
might seem alike to some of the results in the
preceding paper, the truth is that they are much more powerful and the proofs developed are completely new.
Without a doubt, the main contribution of this paper is Theorem \ref{t:caractren}, that provides a
characterization of the hyperbolicity of trains in terms of the
behavior of a real function with two integer parameters.
(In \cite{APRT} we give either necessary or sufficient conditions, but there are no characterizations). This
function describes somehow the distances between some remarkable
geodesics (called \emph{fundamental geodesics}) in the train.
At first sight, Theorem \ref{t:caractren} might not seem very user-friendly.
However, in practice, this tool let us deduce a result about stability of hyperbolicity, even for cases when the
original surface and the modified one are not quasi-isometric (see
Theorem \ref{t:HypStab}).

Theorem \ref{t:caractren} also allows to deduce
both sufficient and necessary conditions
that either guarantee or discard hyperbolicity
(see Theorems \ref{t:suficiente}, \ref{t:necesario1} and \ref{t:necesario2}).
Besides, these three theorems give a much simpler characterization
than Theorem \ref{t:caractren}
for an interesting case of trains:
those for which the lengths of their fundamental geodesics are a quasi-increasing sequence.
We are talking about Theorem \ref{t:caracsel}, another crucial result in this paper.

Theorem \ref{t:union} gives some answers to the following question:
how do some perturbations affect on the hyperbolicity of a flute surface?

For the sake of clarity and readability, we have opted for moving all the technical lemmas to the last section of the paper. This makes the proof of Theorem \ref{t:caractren}, our main result, much more understandable.

\mpb

\noindent
{\bf Notations.}
We denote by $X$ a geodesic metric space.
By $d_X$ and $L_X$ we shall denote,
respectively, the distance and the length in the metric of $X$.
From now on, when there is no possible confusion,
we will not write the subindex $X$.

We denote by $\O$ a train with its Poincar\'e metric.

Given a subset $F$ of the complex plane, we define
$F^+=F \cap \{z\in \CC:\, \Im z\ge 0 \}$, where
$\Im z$ is the imaginary part of $z$.

If $E$ is either a function or a constant related to a domain $\O$,
we will denote by $E'$ or $E^j$ the same function or constant
related to a domain $\O'$ or $\O^j$, respectively.

Finally, we denote
by $c$ and $c_i$, positive constants which can assume
different values in different theorems.

\mpb

\noindent
{\bf Acknowledgements.} We would like to thank Professor J. L. Fern\'andez
for some useful discussions.

\bigskip

\section{Background in Gromov spaces and Riemann surfaces.}

In our study of hyperbolic Gromov spaces we use the notations of \cite{GH}.
We give now the basic facts about these spaces.
We refer to \cite{GH} for more background and further results.

\begin{definicion}
\label{def:1}
Let us fix a point $w$ in a metric space $(X,d)$. We define the
\emph{Gromov product} of $x,y\in X$ with respect to the point $w$
as
$$
(x|y)_w:=\frac12\,\big( d(x,w)+d(y,w)-d(x,y) \big)\ge 0\,.
$$
We say that the metric space $(X,d)$ is $\d$-\emph{hyperbolic}
$(\d\ge 0)$ if
$$
(x|z)_w\ge\min\big\{ (x|y)_w, (y|z)_w \big\}-\d\,,
$$
for every $x,y,z,w\in X$. We say that $X$ is \emph{hyperbolic}
(in the Gromov sense) if the value of $\d$ is not important.
\end{definicion}

It is convenient to remark that this definition of hyperbolicity
is not universally accepted, since sometimes the word hyperbolic
refers to negative curvature or to the existence of Green's
function. However, in this paper we only use the word
\emph{hyperbolic} in the sense of Definition \ref{def:1}.

\spb

\noindent
{\bf Examples:}

\begin {enumerate}
\item Every bounded metric space $X$ is $(diamX)$-hyperbolic (see e.g. \cite[p. 29]{GH}).

\item Every complete simply connected Riemannian manifold with
sectional curvature which is bounded from above by $-k$, with
$k>0$, is hyperbolic (see e.g. \cite[p. 52]{GH}).

\item Every tree with edges of arbitrary length is $0$-hyperbolic
(see e.g. \cite[p. 29]{GH}).
\end{enumerate}

\begin{definicion}
\label{def:geo}
If $\g:[a,b]\longrightarrow X$ is a continuous curve in a metric space $(X,d)$,
the \emph{length} of $\g$ is
$$
L(\g):=\sup\Big\{ \sum_{i=1}^n d(\g(t_{i-1}),\g(t_{i})):\,
a=t_0<t_1<\cdots <t_n=b\Big\}\,.
$$
We say that $\g$ is a \emph{geodesic} if it is an isometry, i.e.
$L(\g|_{[t,s]})=d(\g(t),\g(s))=|t-s|$ for every $s,t\in [a,b]$. We
say that $X$ is a \emph{geodesic metric space} if for every $x,y\in
X$ there exists a geodesic joining $x$ and $y$; we denote by $[x,y]$
any of such geodesics (since we do not require uniqueness of
geodesics, this notation is ambiguous, but convenient as well).
\end{definicion}

\begin{definicion}
\label{def:Rips}
Consider a geodesic metric space $X$.
If $x_1,x_2,x_3\in X$, a
\emph{geodesic triangle} $T=\{x_1,x_2,x_3\}$ is the union of three
geodesics $[x_1,x_2]$, $[x_2,x_3]$ and $[x_3,x_1]$.
We say that $T$ is
$\d$-\emph{thin} if for every $x\in [x_i,x_j]$ we have that
$d(x, [x_j,x_k] \cup [x_k,x_i])\le \d$.
The space $X$
is $\d$-\emph{thin} $($or satisfies the \emph{Rips condition}
with constant $\d)$ if every geodesic triangle in $X$ is
$\d$-thin.
\end{definicion}

As the following basic result states, hyperbolicity is equivalent to
Rips condition:

\begin{teorema}
\label{th:A}
(\cite[p. 41]{GH}) Let us consider a geodesic metric space $X$.

\begin{enumerate}
\item[$(1)$] If $X$ is $\d$-hyperbolic, then it is $4\d$-thin.

\item[$(2)$] If $X$ is $\d$-thin, then it is $4\d$-hyperbolic.
\end{enumerate}
\end{teorema}

\spb

A {\it non-exceptional} Riemann surface $S$ is a
Riemann surface whose universal covering space is
the unit disk $\DD=\{z\in\CC:\; |z|<1\}$, endowed with its Poincar\'e
metric, i.e. the
metric obtained by projecting the Poincar\'e metric of the unit disk
$ds =2 |dz|/(1-|z|^2)$.
Therefore, any simply connected subset of $S$ is isometric to a
subset of $\DD$. With this metric, $S$ is a geodesically complete
Riemannian manifold with constant curvature $-1$, and therefore $S$
is a geodesic metric space. The only Riemann surfaces which are left
out are the {\emph{exceptional} Riemann surfaces, that is to say,
the sphere, the plane, the punctured plane and the tori. It is easy
to study the hyperbolicity of these particular cases.
The Poincar\'e metric is natural and useful in Complex Analysis:
for instance, any holomorphic function between two domains
is Lipschitz with constant $1$,
when we consider the respective Poincar\'e metrics.

\spb

A \emph{Denjoy domain} is a domain $\O$ in the Riemann sphere
with $\p \O\subset \RR\cup\{\infty\}$.
As we mentioned in the introduction of this paper, Denjoy domains
are becoming more and more interesting in Geometric Function Theory
(see e.g. \cite{A}, \cite{APR}, \cite{GJ}, \cite{G}).

\spb

It is obvious that as we focus on more particular kind of surfaces,
we can obtain more powerful results. For this reason we
introduce now a new type of space.

\spb

We have used the word {\it geodesic} in the sense of Definition \ref{def:geo},
that is to say, as a global geodesic or a minimizing geodesic;
however, we need now to deal with a special type of local geodesics:
simple closed geodesics, which obviously can not be minimizing geodesics.
We will continue using the word geodesic with the meaning of Definition \ref{def:geo},
unless we are dealing with closed geodesics.

\begin{definicion}
\label{def:Train}
A \emph{train} is a Denjoy domain $\O \subset \CC$ with $\O \cap
\RR = \cup_{n=0}^{\infty} (a_n,b_n)$, such that $-\infty \le a_0$
and $b_n \le a_{n+1}$ for every $n$.
A \emph{flute surface} is a train with $b_n = a_{n+1}$ for every $n$.

We say that a curve in a train $\O$ is a \emph{fundamental
geodesic} if it is a simple closed geodesic which just intersects
$\RR$ in $(a_0,b_0)$ and $(a_n,b_n)$ for some $n>0$; we denote by
$\g_n$ the fundamental geodesic corresponding to $n$ and
$2l_n:=L_{\O}(\g_n)$. A curve in a train $\O$ is a \emph{second
fundamental geodesic} if it is a simple closed geodesic which
just intersects $\RR$ in $(a_n,b_n)$ and $(a_{n+1},b_{n+1})$ for
some $n \ge 0$; we denote by $\s_n$ the second fundamental
geodesic corresponding to $n$ and $2r_n:=L_{\O}(\s_n)$ (see figure below). If $b_n =
a_{n+1}$, we define $\s_n$ as the puncture at this point and $r_n=0$.
Given $z\in \O$, we define the
\emph{height} of $z$ as $h(z):=d_\O(z,(a_0,b_0))$.
\end{definicion}

\begin{figure}[h]
\label{fig:geodfund}
    \centering
    \psfrag{a0}{$a_0$} %
    \psfrag{-inf}{$-\infty$} %
    \psfrag{b0}{$b_0$} %
    \psfrag{a1}{$a_1$} %
    \psfrag{b1}{$b_1$} %
    \psfrag{a2}{$a_2$} %
    \psfrag{b2}{$b_2$} %
    \psfrag{a3}{$a_3$} %
    \psfrag{b3}{$b_3$} %
    \psfrag{a4}{$a_4$} %
    \psfrag{b4}{$b_4$} %
    \psfrag{g2}{$\g_2$} %
    \psfrag{g3}{$\g_3$} %
    \psfrag{s3}{$\s_3$} %
    \subfigure[Train seen as a subset of the complex plane.]{\epsfig{figure=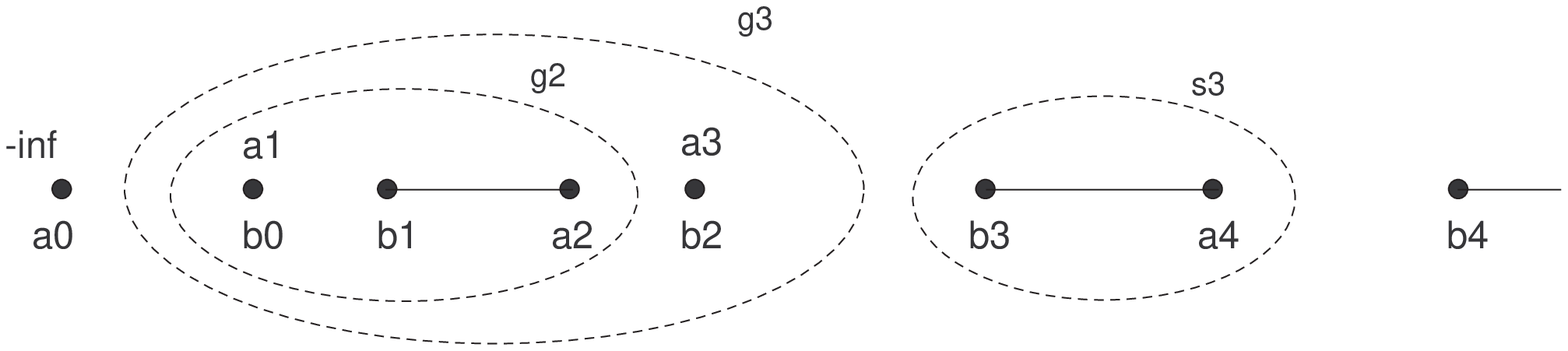,width=10cm}} \\
    \vspace{1cm}
    \subfigure[The same train seen with ``Euclidean
eyes".]{\epsfig{figure=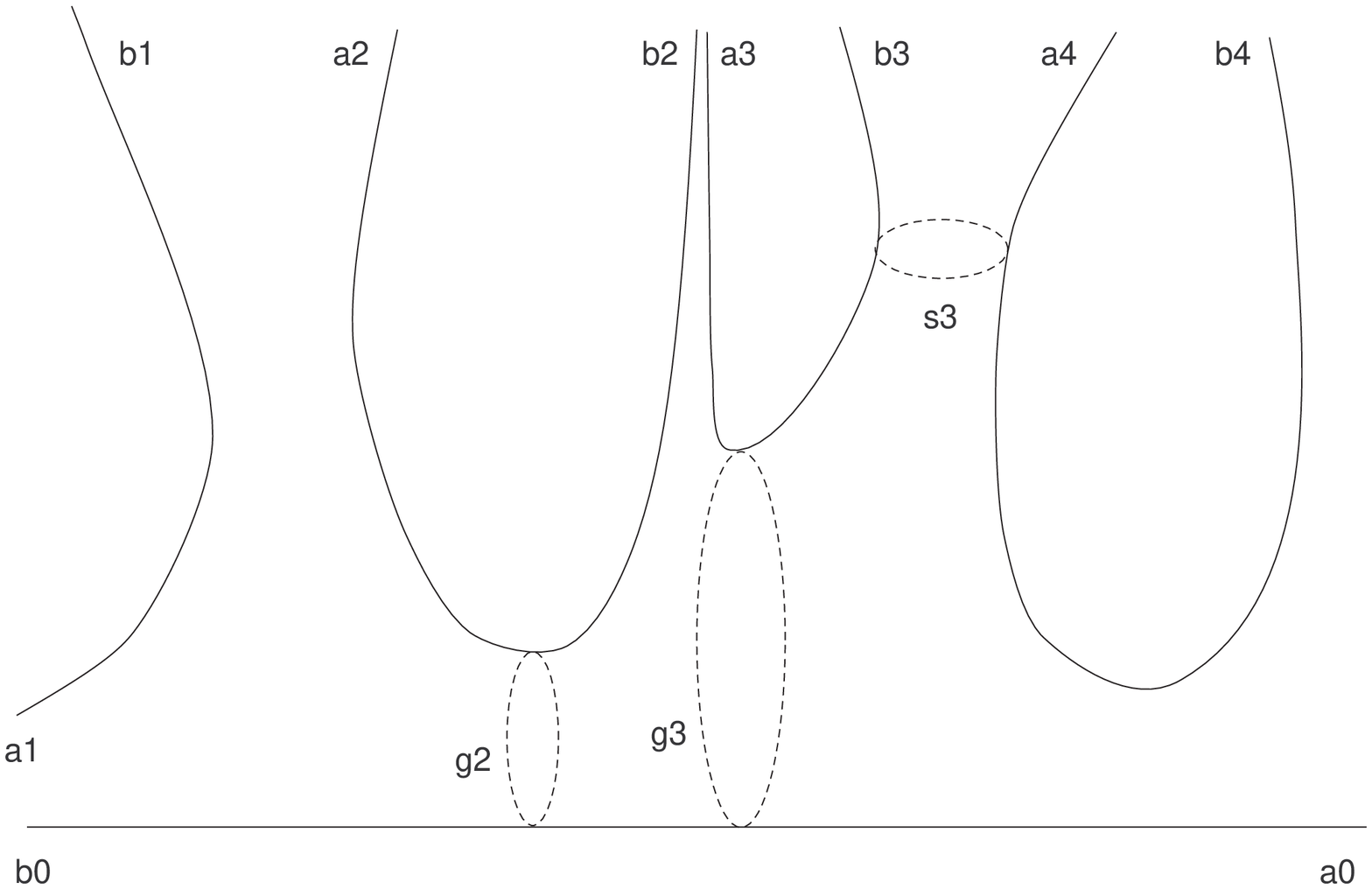,height=5cm, width=10cm}}%
\end{figure}

\noindent
{\bf Remark.}
Recall that in every free homotopy class there exists a single simple closed geodesic, assuming that punctures are simple closed geodesics with length equal to zero. That is why both the fundamental geodesic and the second fundamental geodesic are unique for every $n$.

A train is a flute surface if and only if
every second fundamental geodesic is a puncture.

Flute surfaces are the simplest examples of
infinite ends; furthermore, in a flute surface it is possible to give a fairly
precise description of the ending geometry (see, e.g. \cite{H}).

\bigskip

\section{The main results.}
\label{s:MainResults}

\medskip

It is not difficult to see that the values of
$\{l_n\}$ and $\{r_n\}$ determine a train, since for every $n$ there exists a single fundamental geodesic and a single second fundamental geodesic (see the Remark to Definition \ref{def:Train}).
Then, there must exist a characterization of hyperbolicity
in terms of the lengths of the fundamental geodesics. It would be desirable to obtain such a characterization, since these lengths describe the Denjoy domain from a simple geometric viewpoint.

In order to obtain this characterization, we need to introduce the following
functions.

(We refer to the next section for the details of the proofs of technical lemmas. We think that this structure makes the paper more readable, because it shortens considerably the proof of Theorem \ref{t:caractren}).

\begin{definicion}
Let us consider a sequence of positive numbers
$\{l_n\}_{n=1}^\infty$
and a sequence of non-negative numbers
$\{r_n\}_{n=1}^\infty$.
Consider $n\ge 1$ and $0\le h \le l_n$.
We define $A_n(h):=\max\{m<n: \, l_m\le h \}$ if this set is non-empty
and $A_n(h):=1$ in other case,
$B_n(h):=\min\{m>n: \, l_m\le h \}$ if this set is non-empty
and $B_n(h):=\infty$ in other case,
$$
\D(k):= e^{-l_{k}} + e^{-l_{k+1}} +
e^{-\frac12
(l_{k}+l_{k+1}-r_k)_+} + (r_k-l_{k}-l_{k+1})_+ \,,
$$
and
$$
\Gamma_{nm}(h)
:=
\left\{ \begin{array}{ll}
\!\!
\big(r_{m} + h - l_{m+1} \big)_+  +
e^{h} \!\!\! \displaystyle\sum_{k=m+1}^{n-1} \!\!\! \D(k) \,,
& \;
\; \text{ if } \, m<n\, \text{ and } \, l_m \le h \,,
\\
l_m - h +
e^{h} \displaystyle\sum_{k=m}^{n-1} \D(k) \,,
& \;
\; \text{ if } \, m<n \, \text{ and } \, l_m> h \,,
\\
\, & \,
\\
\min  \big\{h,\, l_n - h \big\} \,,
& \;
\; \text{ if } \, m=n \,,
\\
\, & \,
\\
l_m - h +
e^{h} \displaystyle\sum_{k=n}^{m-1} \D(k) \,,
& \;
\; \text{ if } \, m>n \, \text{ and } \, l_m> h \,,
\\
\big(r_{m-1} + h - l_{m-1} \big)_+  +
e^{h} \displaystyle\sum_{k=n}^{m-2}  \D(k) \,,
& \;
\; \text{ if } \, m>n \, \text{ and } \, l_m \le h \,.
\end{array}
\right.
$$
The functions $\Gamma_{nm}(h)$ are naturally associated to trains
by taking $\{l_n\}_{n=1}^\infty$ and $\{r_n\}_{n=1}^\infty$
as the half-lengths of their fundamental geodesics.
\end{definicion}

\begin{teorema}
\label{t:caractren}
A train $\O$ is hyperbolic if and only if
$$
K:=\sup_{n\ge 1} \sup_{h\in [0,l_n]} \min_{m\in [A_n(h),B_n(h)]}
\Gamma_{nm}(h) < \infty\,.
$$
Furthermore,
if $\O$ is $\d$-hyperbolic, then $K$ is bounded by a constant
which only depends on $\d$;
if $K<\infty$, then $\O$ is $\d$-hyperbolic, with $\d$ a constant
which only depends on $K$.
\end{teorema}

\noindent
{\bf Remarks.}

$(1)$
Notice that this is a real variable characterization of the
hyperbolicity, although the hyperbolicity is a concept of complex
geometry, since we consider the Poincar\'e metric in each train.

$(2)$
Theorem \ref{t:caractren}
clearly improves \cite[Theorem 5.3]{APRT}:
we need to know the lengths of the fundamental geodesics instead of
the precise location of these geodesics and the distances to
$\RR$ from their points.

$(3)$
The proof of Theorem \ref{t:caractren}
gives that its conclusion also holds if we replace $K$ by
$$
K(l_0):=\sup_{n\ge 1} \sup_{h\in [l_0,l_n]} \min_{m\in [A_n(h),B_n(h)]}
\Gamma_{nm}(h) < \infty\,,
$$
for any fixed $l_0>0$. In this case, the constant $\d$
depends on $K(l_0)$ and $l_0$.

\begin{proof}
By \cite[Theorem 5.3]{APRT}, $\O$ is $\d$-hyperbolic if and only if
$$
K_1:=\sup_{n\ge 1} \sup_{z\in \g_n} \inf_{m\ge 0}
d_{\O} \big(z,(a_m,b_m)\big) < \infty\,,
$$
with the appropriate dependence of the constants
(if $\O$ is $\d$-hyperbolic, then $K_1$ is bounded by a constant
which only depends on $\d$;
if $K_1<\infty$, then $\O$ is $\d$-hyperbolic, with $\d$ a constant
which only depends on $K_1$).

Fix any constant $l_0>0$.
Notice that:

\smallskip

$(1)$ $d_{\O} \big(z,(a_0,b_0)\big) = h(z)$
and $d_{\O} \big(z,(a_n,b_n)\big) = l_n - h(z)$.
Since any $z$ with $h(z)<l_0$ verifies
$$
\inf_{m\ge 0} d_{\O} \big(z,(a_m,b_m)\big)
\le d_{\O} \big(z,(a_0,b_0)\big) = h(z) < l_0\,,
$$
we only need to consider $z$ with $l_0 \le h(z) \le l_n$.

\smallskip

From now on, let us fix $n\ge 1$ and $z\in \g_n$ with $l_0 \le h(z) \le l_n$.

\smallskip

$(2)$ If $k<m<n$, with $l_m\le h(z)$, let us consider
the geodesic $\sigma$ which gives the minimum distance between
$z$ and $(a_k,b_k)$.
Define the point $w:=\sigma \cap \g_m$; hence
$d_{\O} \big(z,w\big) <  d_{\O} \big(z,(a_k,b_k)\big)$ and
Lemma \ref{l:cartesiana} gives
$$
d_{\O} \big(z,(a_m,b_m)\big)
\le d_2 \big(z,(a_m,b_m) \cap \g_m\big)
\le d_2 \big(z,w \big)
\le 3 \, d_{\O} \big(z,w \big)
< 3 \, d_{\O} \big(z,(a_k,b_k)\big)\,.
$$
In a similar way, if $k>m>n$, with $l_m\le h(z)$, then
$d_{\O} \big(z,(a_m,b_m)\big) < 3 \,d_{\O} \big(z,(a_k,b_k)\big)$.
Hence we only need to consider $d_{\O} \big(z,(a_m,b_m)\big)$
with $m\in \{0\}\cup [A_n(h(z)),B_n(h(z))]$,
in order to study if $K_1$ is finite.

\smallskip

$(3)$ If $m\in (A_n(h(z)),n)$, then $l_0 \le h(z) < l_m$.
By Lemma \ref{l:cartesiana2},
we can replace
$d_{\O} \big(z,(a_m,b_m)\big)$ by
$d_1 \big(z,\g_m \cap (a_m,b_m)\big)$.
If $z_m$ is the point in $\g_m$ with $h(z_m)=h(z)$, then
$d_1 \big(z,\g_m \cap (a_m,b_m)\big):= d_{\O} (z,z_m) + l_m - h(z)$.
Standard hyperbolic trigonometry in quadrilaterals (see e.g. \cite[p. 88]{F})
gives that
$$
 d_\O(z,z_m) = 2 \Arcsinh \Big( \sinh \frac12\, d_\O(\g_m,\g_n) \cosh h(z) \Big) \,.
$$
Recall that $(a_0,b_0)$ contains the shortest geodesic joining $\g_m$ and $\g_n$.
By Corollary \ref{c:aex} we can replace $d_\O (z,z_m)$
by $d_\O(\g_m,\g_n) \, e^{h(z)}$,
and therefore $d_1 \big(z,\g_m \cap (a_m,b_m)\big)$
by $d_\O(\g_m,\g_n) \, e^{h(z)} + l_m - h(z)$.
Standard hyperbolic trigonometry in right-angled hexagons (see e.g. \cite[p. 86]{F})
gives that
$$
d_\O(\g_k,\g_{k+1})
= \Arccosh \frac{\cosh r_k + \cosh l_{k} \cosh l_{k+1}}{\sinh l_{k} \sinh l_{k+1}}
$$
for every $k \ge 1$.
Proposition \ref{p:f} gives
$$
d_\O(\g_k,\g_{k+1})= f(l_{k},l_{k+1},r_{k})
\asymp  e^{-l_{k}} + e^{-l_{k+1}} + e^{-\frac12 (l_{k}+l_{k+1}-r_k)_+}
+ (r_k-l_{k}-l_{k+1})_+ = \D(k) \,,
$$
for every $k \in (A_n(h(z)),n)$, since then $l_{k},l_{k+1} \ge h(z) \ge l_0$.
Therefore we can replace
$d_{\O} \big(z,(a_m,b_m)\big)$
by
$$
l_m - h(z) +
e^{h(z)} \sum_{k=m}^{n-1} \D(k) \,.
$$
A symmetric argument gives that if $m\in (n,B_n(h(z)))$,
then we can replace $d_{\O} \big(z,(a_m,b_m)\big)$ by
$$
l_m - h(z) +
e^{h(z)} \sum_{k=n}^{m-1} \D(k) \,.
$$
\indent
$(4)$ If $m=A_n(h(z))$, then $h(z) \ge l_m$.
If $z_{m+1}$ is the point in $\g_{m+1}$ with $h(z_{m+1})=h(z)$,
by Lemma \ref{l:partirendos}, we can replace
$d_{\O} \big(z,(a_m,b_m)\big)$ by
$d_{\O} \big(z,z_{m+1}\big) + d_{\O} \big(z_{m+1},(a_m,b_m)\big)$.
We have seen in $(3)$ that
we can replace $d_{\O} \big(z,z_{m+1}\big)$ by
$$
e^{h(z)}  \!\! \sum_{k=m+1}^{n-1} \!\!  \D(k) \,.
$$
Standard hyperbolic trigonometry in pentagons (see e.g. \cite[p. 87]{F})
gives that
$$
\sinh d_{\O} \big(z_{m+1},(a_m,b_m)\big)
= - \cosh l_{m} \sinh h(z) + \sinh l_{m} \cosh h(z) \cosh d_\O(\g_m,\g_{m+1}) \,.
$$
Standard hyperbolic trigonometry in right-angled hexagons (see e.g. \cite[p. 86]{F})
gives that
$$
\cosh d_\O(\g_m,\g_{m+1})
= \frac{\cosh r_m + \cosh l_{m} \cosh l_{m+1}}{\sinh l_{m} \sinh l_{m+1}} \;,
$$
and hence
$$
\aligned
\sinh d_{\O} \big(z_{m+1},(a_m,b_m)\big)
& = - \cosh l_{m} \sinh h(z) + \cosh h(z) \,
\frac{\cosh r_m + \cosh l_{m} \cosh l_{m+1}}{\sinh l_{m+1}}
\\
& = \frac{
\cosh l_{m} \big(\cosh l_{m+1} \cosh h(z)  - \sinh l_{m+1} \sinh h(z)\big)
+ \cosh r_m \cosh h(z)}{\sinh l_{m+1}}
\\
& = \frac{
\cosh l_{m} \cosh \big( l_{m+1} - h(z) \big) + \cosh r_m \cosh h(z)}{\sinh l_{m+1}}
= \sinh F \big(l_{m}, l_{m+1}, r_m,  h(z) \big) \,,
\endaligned
$$
where $F$ is the function in Proposition \ref{p:F}. Therefore,
Corollary \ref{c:xyth} gives that we can replace $d_{\O}
\big(z_{m+1},(a_m,b_m)\big)$ by $\big(r_m + h(z) - l_{m+1} \big)_+$.
Consequently, we can substitute $d_{\O} \big(z,(a_m,b_m)\big)$ by
$$
\big(r_{m} + h(z) - l_{m+1} \big)_+  +
e^{h(z)} \!\! \sum_{k=m+1}^{n-1} \!\!  \D(k) \,.
$$
\indent
A symmetric argument gives that if $m=B_n(h(z))$,
then we can replace $d_{\O} \big(z,(a_m,b_m)\big)$ by
$$
\big(r_{m-1} + h(z) - l_{m-1} \big)_+  +
e^{h(z)}  \sum_{k=n}^{m-2}  \D(k) \,.
$$

Notice that each time that we replace a quantity by another in this proof,
the constants are under control.
Let us remark that
$(1)$, $(2)$, $(3)$ and $(4)$ give the result, with
$\inf_{m\in [A_n(h),B_n(h)]}\Gamma_{nm}(h)$ instead of
$\min_{m\in [A_n(h),B_n(h)]}\Gamma_{nm}(h)$.

Let us see now that this infimum is attained. Seeking for a
contradiction, suppose that the latest statement is not true.
Therefore, $B_n(h) = \infty$ and $l_m
> h$ for every $m > n$. Then, there exists an increasing sequence of
integer numbers $\{m_j\}$ with $\lim_{j\to \infty}\Gamma_{nm_j}(h) =
\inf_{m\in [A_n(h),\infty)}\Gamma_{nm}(h)$. By choosing a
subsequence if it is necessary, we can assume that
$\{\Gamma_{nm_j}(h)\}_j$ is a decreasing sequence. Hence,
$$
\Gamma_{nm_{j+1}}(h)
= l_{m_{j+1}} - h +  e^{h} \displaystyle\sum_{k=n}^{m_{j+1}-1} \D(k)
< \Gamma_{nm_j}(h)
= l_{m_j} - h +  e^{h} \displaystyle\sum_{k=n}^{m_j-1} \D(k) \,.
$$
Consequently, we have that $l_{m_{j+1}} < l_{m_{j}} < l_{m_{1}}$ for every $j$,
and
$$
\Gamma_{nm_j}(h)
 = l_{m_j} - h +  e^{h} \displaystyle\sum_{k=n}^{m_j-1} \D(k)
\ge e^{h} \displaystyle\sum_{k=n}^{m_j}  e^{-l_{k}}
\ge e^{h} \displaystyle\sum_{k=1}^{j}  e^{-l_{m_k}}
\ge e^{h} \, j\,  e^{-l_{m_1}}  .
$$
Hence, $\lim_{j\to \infty}\Gamma_{nm_j}(h)
= \lim_{j\to \infty} e^{h} \, j\, e^{-l_{m_1}} = \infty$,
which is a contradiction.
This finishes the proof.
\end{proof}

\begin{lema}
\label{l:crec1}
For every $r_{k}\ge 0$ and $0 < l_{k} \le h \le l_{k+1}$, we have
$$
\big(r_{k} + h - l_{k+1} \big)_+
< e^{h}  \D(k) \,.
$$
\end{lema}

\begin{proof}
Let us remark that it is sufficient to prove
$$
r_{k} + h - l_{k+1}
< e^{h}  \big( e^{-\frac12
(l_{k}+l_{k+1}-r_k)_+} + (r_k-l_{k}-l_{k+1})_+ \big) \,,
$$
for every $r_{k}\ge 0$ and $0 < l_{k} \le h \le l_{k+1}$.

Since the left hand side of the inequality does not depend on $l_k$
and the right hand side is a decreasing function on $l_k$, it is
sufficient to prove
$$
r_{k} + h - l_{k+1}
< e^{h}  \big( e^{-\frac12
(h+l_{k+1}-r_k)_+} + (r_k-h-l_{k+1})_+ \big) \,,
$$
for every $r_{k}\ge 0$ and $0 < h \le l_{k+1}$.

If $r_k\le h+l_{k+1}$, then the inequality is
$$
r_{k} + h - l_{k+1}
< e^{h}  e^{-\frac12 (h+l_{k+1}-r_k)}
= e^{\frac12 (r_k+h-l_{k+1})} \,,
$$
which trivially holds since $t< e^{t/2}$ for every real number $t$.

If $r_k\ge h+l_{k+1}$, then the inequality is
$$
r_{k} + h - l_{k+1}
< e^{h} (1+r_k-h-l_{k+1}) \,.
$$
Since $h>1$, it is clear that the function
$$
U(r_{k}):= e^{h} (1+r_k-h-l_{k+1})
- r_{k} - h + l_{k+1}
$$
is increasing in $r_k\in [h + l_{k+1}, \infty)$.
Then $U(r_{k}) \ge U(h+l_{k+1}) = e^h -2h > 0$,
and the inequality holds.
\end{proof}

\begin{prop}
\label{p:minimo}
In any train $\O$ we have
$$
\min_{m\in [A_n(h),B_n(h)]} \Gamma_{nm}(h)
= \min_{m\ge 1} \Gamma_{nm}(h) \,,
$$
for every $n\ge 1$ and $0\le h \le l_n$.
\end{prop}

\begin{proof}
Fix $n\ge 1$ and $0\le h \le l_n$.
If $m< A_n(h)$, then Lemma \ref{l:crec1} gives
$\Gamma_{nm}(h) > \Gamma_{nA_n(h)}(h)$:
$$
\aligned
\Gamma_{nm}(h)
&
\ge e^h \!\! \sum_{k=m+1}^{n-1} \!\! \D(k)
\ge e^h \!\!\!\! \sum_{k=A_n(h)}^{n-1} \!\!\!\! \D(k)
= e^h \D(A_n(h)) + e^h \!\! \!\! \!\! \sum_{k=A_n(h)+1}^{n-1} \!\!\!\! \!\!  \D(k)
\\
&
> \big(r_{A_n(h)} + h - l_{A_n(h)+1} \big)_+
+ e^h \!\! \!\!\!\! \sum_{k=A_n(h)+1}^{n-1} \!\!\!\! \!\! \D(k)
= \Gamma_{nA_n(h)}(h) \,.
\endaligned
$$
The case $m > B_n(h)$ is similar.
\end{proof}

\begin{prop}
\label{p:altura}
If for some $n$ we have $l_m\ge l_n$ for every $m \ge n$,
then the conclusion of Theorem $\ref{t:caractren}$
also holds if we replace $[A_n(h),B_n(h)]$ by $[A_n(h),n]$ for this $n$.
\end{prop}

\begin{proof}
It suffices to remark that
for every $z\in \g_n$ and $m > n$, we have
$d_\O (z, (a_n,b_n)) = l_n- h(z) \le l_m- h(z) < d_\O (z, (a_m,b_m))$.
\end{proof}

Although to compute the minimum and the supremum
in Theorem \ref{t:caractren} can be difficult in the general case,
Theorem \ref{t:caractren} is the main tool in order to obtain the remaining
results of this paper.
We start with an elementary corollary.

\begin{prop}
\label{p:lacotado}
Let us consider a train $\O$ with $l_n\le c$ for every $n$.
Then $\O$ is $\d$-hyperbolic, where $\d$ is a constant which
only depends on $c$.
\end{prop}

\begin{proof}
For each positive integer $n$, we have
$\Gamma_{nn}(h) := \min  \big\{h,\, l_n - h \big\} \le l_n \le c$
for every $h \in [0,l_n]$.
Hence, $K \le c$ and Theorem \ref{t:caractren} finishes the proof.
\end{proof}

One of the important problems in the study of any property
is to obtain its stability under appropriate deformations.
Theorem \ref{t:caractren} allows to prove a result which shows that
hyperbolicity is stable under bounded perturbations of the lengths
of the fundamental geodesics. Theorem \ref{t:HypStab} is
particularly remarkable since there are very few results on
hyperbolic stability which do not involve quasi-isometries. We need
a previous lemma; it deals with some kind of reverse inequality to
the one in Lemma \ref{l:crec1}.

\begin{lema}
\label{l:crec2}
For every $r_{k},l_{k+1}\ge 0$ and $0 \le h \le l_{k}$, we have
$$
e^{h} \big( e^{-\frac12 (l_{k}+l_{k+1}-r_k)_+} + (r_k-l_{k}-l_{k+1})_+ \big)
\le \big( 1 + (r_{k} + h - l_{k+1} )_+ \big)
\, e^{\frac12 (r_{k} + h - l_{k+1} )_+} .
$$
\end{lema}

\begin{proof}
Since the right hand side of the inequality does not depend on $l_k$
and the left hand side is a decreasing function on $l_k$, it is
sufficient to prove
$$
e^{h} \big( e^{-\frac12 (h+l_{k+1}-r_k)_+} + (r_k-h-l_{k+1})_+ \big)
\le ( 1 + \big(r_{k} + h - l_{k+1} )_+ \big)
\, e^{\frac12 (r_{k} + h - l_{k+1} )_+} .
$$
for every $r_{k},l_{k+1}, h\ge 0$.

If $h+l_{k+1}-r_k\ge 0$, the inequality is direct since
$$
e^{h} \big( e^{-\frac12 (h+l_{k+1}-r_k)_+} + (r_k-h-l_{k+1})_+ \big)
= e^{h} e^{-\frac12 (h+l_{k+1}-r_k)}
= e^{\frac12 (r_{k} + h - l_{k+1} )} .
$$
\indent
If $h+l_{k+1}-r_k < 0$, then $r_k-l_{k+1} > h$
and $(r_{k} + h - l_{k+1})_+ > 2h$; consequently,
$$
e^{h} \big( e^{-\frac12 (h+l_{k+1}-r_k)_+} + (r_k-h-l_{k+1})_+ \big)
= e^{h} \big( 1 + r_k-h-l_{k+1} \big)
< ( 1 + \big(r_{k} + h - l_{k+1} )_+ \big)
\, e^{\frac12 (r_{k} + h - l_{k+1} )_+} .
$$
\end{proof}

Next, the result about stability that we have talked about before
Lemma \ref{l:crec2}.
Theorem \ref{t:HypStab} is both
a qualitative and a quantitative result.

\begin{teorema}
\label{t:HypStab}
Let us consider two trains $\O$, $\O'$ and a constant $c$
such that $|r'_n-r_n| \le c$, and $|l'_{n}-l_{n}|\le c$ for every $n\ge 1$.
Then $\O$ is hyperbolic if and only if $\O'$ is hyperbolic.

\noindent
Furthermore, if $\O$ is $\d$-hyperbolic, then $\O'$ is $\d'$-hyperbolic,
with $\d'$ a constant which only depends on $\d$ and $c$.
\end{teorema}

This result is a significant improvement with respect to \cite[Theorem 5.33]{APRT}, since, in that paper, the lengths $r_n$ and $r_n'$ were required to be bounded, whereas Theorem \ref{t:HypStab} only requires  $r_n-r_n'$ to be bounded. Notice that this is a much weaker condition. Furthermore, the argument in the proof is completely new.

\medskip

\noindent
{\bf Remarks.}

$(1)$
Notice that in many cases $\O$ and $\O'$ are not quasi-isometric
(for example, if there exists a subsequence $\{n_k\}_k$ with
$\lim_{k\to\infty}l_{n_k}=0$ and $l_{n_k}'\ge c_0>0$).

$(2)$
We have examples which show that Theorem \ref{t:HypStab} is sharp:
if we change the constants in Theorem \ref{t:HypStab} by any function
growing slowly to infinity, then the conclusion of Theorem \ref{t:HypStab}
does not hold.
For instance,
if $\{r_n\}$ is bounded and $\{r_n'\}$ is not bounded,
then there exists $\{l_n\}=\{l_n'\}$
with $\O$ hyperbolic and $\O'$ not hyperbolic.

\begin{proof}
By symmetry, it is sufficient to prove that if $\O$ is $\d$-hyperbolic,
then $\O'$ is $\d'$-hyperbolic,
with $\d'$ a constant which only depends on $\d$ and $c$.
Therefore, let us assume that $\O$ is $\d$-hyperbolic.

Notice that
$e^{-l_{k}} + e^{-l_{k+1}} \le e^c \big( e^{-l'_{k}} + e^{-l'_{k+1}} \big)$.

If $l_{k}+l_{k+1} \le r_k$, then
$e^{-\frac12 (l_{k}+l_{k+1}-r_k)_+} + (r_k-l_{k}-l_{k+1})_+
=1 + r_k-l_{k}-l_{k+1}$
and
$$
e^{-\frac12 (l'_{k}+l'_{k+1}-r'_k)_+} + (r'_k-l'_{k}-l'_{k+1})_+
\le 1 + 3c + r_k-l_{k}-l_{k+1}
\le (1 + 3c) \big( e^{-\frac12 (l_{k}+l_{k+1}-r_k)_+} + (r_k-l_{k}-l_{k+1})_+ \big) \,.
$$
\indent
If $l'_{k}+l'_{k+1} \ge r'_k$, then
$$
e^{-\frac12 (l'_{k}+l'_{k+1}-r'_k)_+} + (r'_k-l'_{k}-l'_{k+1})_+
= e^{-\frac12 (l'_{k}+l'_{k+1}-r'_k)_+}
\le e^{3c/2} \big( e^{-\frac12 (l_{k}+l_{k+1}-r_k)_+} + (r_k-l_{k}-l_{k+1})_+ \big) \,.
$$
\indent
If $l_{k}+l_{k+1} > r_k$ and $l'_{k}+l'_{k+1} < r'_k$, then
$$
\aligned
l_{k}+l_{k+1} - r_k & \le l'_{k}+l'_{k+1} - r'_k +3c < 3c \,,
\\
r'_k-l'_{k}-l'_{k+1} & \le r_k-l_{k}-l_{k+1} +3c < 3c \,,
\endaligned
$$
and consequently
$$
\aligned
e^{-\frac12 (l'_{k}+l'_{k+1}-r'_k)_+} + (r'_k-l'_{k}-l'_{k+1})_+
& = 1 + r'_k-l'_{k}-l'_{k+1}
< (1 + 3c) \, e^{3c/2} e^{-3c/2}
\\
&
< (1 + 3c) \, e^{3c/2} \big( e^{-\frac12 (l_{k}+l_{k+1}-r_k)_+}
+ (r_k-l_{k}-l_{k+1})_+ \big) \,.
\endaligned
$$
\indent
Therefore
$$
e^{-l'_{k}} + e^{-l'_{k+1}} + e^{-\frac12 (l'_{k}+l'_{k+1}-r'_k)_+} +
(r'_k-l'_{k}-l'_{k+1})_+
\le (1 + 3c) \, e^{3c/2} \big( e^{-l_{k}} + e^{-l_{k+1}}
+ e^{-\frac12 (l_{k}+l_{k+1}-r_k)_+} + (r_k-l_{k}-l_{k+1})_+ \big) \,,
$$
i.e. $\D'(k)  \le (1 + 3c) \, e^{3c/2} \D(k)$.
We also have
$$
\aligned
(r'_m+h-l'_{m+1})_+ & \le 2c + (r_m+h-l_{m+1})_+ \,,
\\
l'_m-h & \le c + l_m - h \,,
\\
\min  \big\{h,\, l'_n - h \big\} & \le c + \min  \big\{h,\, l_n - h \big\} \,.
\endaligned
$$
Hence, we conclude
$$
\big(\Gamma_{nm}\big)'(h) \le (1 + 3c) \, e^{3c/2} \Gamma_{nm}(h) + 2c \,,
$$
for every $n,m\ge 1$ and $h\ge 0$ with either
$m=n$ or $l_m,l'_m \le h$ or $l_m,l'_m > h$.

We deal now with the other cases.
Let us assume that $m \in [A'_n(h),n)$.
The case $m \in (n,B'_n(h)]$ is similar.

If $l'_m \le h < l_m$, then $m= A'_n(h)$ and $l'_m \le h < l'_{m+1}$.
Applying Lemma \ref{l:crec1} we obtain
$$
\aligned
\big(\Gamma_{nm}\big)'(h)
& = \big(r'_{m} + h - l'_{m+1} \big)_+  +
e^{h} \!\!\! \sum_{k=m+1}^{n-1} \!\!\! \D'(k)
< e^{h} \sum_{k=m}^{n-1} \D'(k)
\\
& \le l_m-h + (1 + 3c) \, e^{3c/2}  e^{h} \sum_{k=m}^{n-1} \D(k)
\le (1 + 3c) \, e^{3c/2} \Gamma_{nm}(h) \,.
\endaligned
$$
\indent
If $l_m \le h < l'_m$, then $m > A'_n(h)$ and $h < l'_{m+1}$.
We also have $l'_m - h \le l'_m - l_m \le c$.
Applying Lemma \ref{l:crec2} we obtain
$$
\aligned
\big(\Gamma_{nm}\big)'(h)
& = l'_{m} - h
+ e^{h-l'_{m}} + e^{h-l'_{m+1}}
+ e^{h} \big( e^{-\frac12 (l'_{m}+l'_{m+1}-r'_m)_+} + (r'_m-l'_{m}-l'_{m+1})_+ \big)
+ e^{h} \!\!\! \sum_{k=m+1}^{n-1} \!\!\! \D'(k)
\\
& \le c + 2
+ \big( 1 + (r'_{m} + h - l'_{m+1} )_+ \big) \, e^{\frac12 (r'_{m} + h - l'_{m+1} )_+}
+ (1 + 3c) \, e^{3c/2}  e^{h} \!\!\! \sum_{k=m+1}^{n-1} \!\!\! \D(k)
\\
& \le c + 2
+ \big( 1 + 2c + (r_{m} + h - l_{m+1} )_+ \big) \, e^c e^{\frac12 (r_{m} + h -
l_{m+1})_+}
+ (1 + 3c) \, e^{3c/2}  e^{h} \!\!\! \sum_{k=m+1}^{n-1} \!\!\! \D(k)
\\
& \le c + 2
+ \big( 1 + 2c + \Gamma_{nm}(h) \big) \, e^c e^{\frac12 \Gamma_{nm}(h)}
+ (1 + 3c) \, e^{3c/2} \Gamma_{nm}(h) \,.
\endaligned
$$
\indent
We can conclude in any case
$$
\aligned
\sup_{h\in [0,\min\{l_n,l_n'\}]} & \min_{m\in [A'_n(h),B'_n(h)]}
\big(\Gamma_{nm}\big)'(h)
 = \sup_{h\in [0,\min\{l_n,l_n'\}]} \min_{m\ge 1}
\big(\Gamma_{nm}\big)'(h)
\\
& \le
\sup_{h\in [0,l_n]} \min_{m\ge 1} \Big( c + 2
+ \big( 1 + 2c + \Gamma_{nm}(h) \big) \, e^c e^{\frac12 \Gamma_{nm}(h)}
+ (1 + 3c) \, e^{3c/2} \Gamma_{nm}(h) \Big)
\\
& \le c + 2
+ \big( 1 + 2c + K \big) \, e^c e^{\frac12 K}
+ (1 + 3c) \, e^{3c/2} K ,
\endaligned
$$
for every $n\ge 1$, where $K$ only depends on $\d$,
by Theorem \ref{t:caractren} and Proposition \ref{p:minimo}.

If for some $n$ we have $l_n < l'_n$ and $h \in [l_n, l'_n]$,
then $\big(\Gamma_{nn}\big)'(h) \le l'_n - h \le l'_n - l_n \le c$ and
$$
\sup_{h\in [l_n,l_n']} \min_{m\in [A'_n(h),B'_n(h)]}
\big(\Gamma_{nm}\big)'(h) \le c \,.
$$
Therefore,
$K' \le c + 2
+ \big( 1 + 2c + K \big) \, e^c e^{\frac12 K}
+ (1 + 3c) \, e^{3c/2} K$,
and the conclusion holds by Theorem \ref{t:caractren}.
\end{proof}

Theorem \ref{t:HypStab} has the following direct consequence.

\begin{corolario}
\label{c:HypStab}
Let us consider two trains $\O$, $\O'$ such that
$r'_n=r_n$, and $l'_{n}=l_{n}$ for every $n\ge N$.
Then $\O$ is hyperbolic if and only if $\O'$ is hyperbolic.
\end{corolario}

Theorems \ref{t:sel0} and \ref{t:sel} are simpler versions of
Theorem \ref{t:caractren}, which can be applied in many occasions,
and are obtained by replacing $\Gamma_{nm}(h)$ for
$\Gamma_{nm}^*(h)$ and $\Gamma_{nm}^0(h)$, respectively. We define
now these functions.

\begin{definicion}
Let us consider a sequence of positive numbers
$\{l_n\}_{n=1}^\infty$
and a sequence of non-negative numbers
$\{r_n\}_{n=1}^\infty$.
Consider $n\ge 1$ and $0\le h \le l_n$.
We define
$$
\Gamma_{nm}^*(h)
:=
\left\{ \begin{array}{ll}
\!\!
\big(r_{m} + h - l_{m+1} \big)_+  +
e^{h} \!\!\! \displaystyle\sum_{k=m+1}^{n} \!\!\! e^{-l_k} \,,
& \;
\; \text{ if } \, m<n\, \text{ and } \, l_m \le h \,,
\\
l_m - h +
e^{h} \displaystyle\sum_{k=m}^{n} e^{-l_k} \,,
& \;
\; \text{ if } \, m<n \, \text{ and } \, l_m> h \,,
\\
\, & \,
\\
\min  \big\{h,\, l_n - h \big\} \,,
& \;
\; \text{ if } \, m=n \,,
\\
\, & \,
\\
l_m - h +
e^{h} \displaystyle\sum_{k=n}^{m} e^{-l_k} \,,
& \;
\; \text{ if } \, m>n \, \text{ and } \, l_m> h \,,
\\
\big(r_{m-1} + h - l_{m-1} \big)_+  +
e^{h} \displaystyle\sum_{k=n}^{m-1}  e^{-l_k} \,,
& \;
\; \text{ if } \, m>n \, \text{ and } \, l_m \le h \,,
\end{array}
\right.
$$
and
$$
\Gamma_{nm}^0(h)
:=
\left\{ \begin{array}{ll}
\!\!
e^{h} \!\!\! \displaystyle\sum_{k=m+1}^{n} \!\!\! e^{-l_k} \,,
& \;
\; \text{ if } \, m<n\, \text{ and } \, l_m \le h \,,
\\
e^{h} \displaystyle\sum_{k=n}^{m-1}  e^{-l_k} \,,
& \;
\; \text{ if } \, m>n \, \text{ and } \, l_m \le h \,,
\\
\Gamma_{nm}^*(h)  \,,
& \;
\; \text{ if } \, m>n \, \text{ in other case. }
\end{array}
\right.
$$
The functions $\Gamma_{nm}^*(h)$ and $\Gamma_{nm}^0(h)$ are naturally associated to
trains by taking $\{l_n\}_{n=1}^\infty$ and $\{r_n\}_{n=1}^\infty$
as the half-lengths of their fundamental geodesics.
\end{definicion}

\begin{teorema}
\label{t:sel0}
Let us consider a train $\O$ such that there exists a constant $c>0$ with
$r_n\le 2c + |l_n-l_{n+1}|$ for every $n \ge 1$.
Then $\O$ is hyperbolic if and only if
$$
K^*:=\sup_{n\ge 1} \sup_{h\in [0,l_n]} \min_{m\in [A_n(h),B_n(h)]}
\Gamma_{nm}^*(h) < \infty\,.
$$
Furthermore,
if $\O$ is $\d$-hyperbolic, then $K^*$ is bounded by a constant
which only depends on $\d$ and $c$;
if $K^*<\infty$, then $\O$ is $\d$-hyperbolic, with $\d$ a constant
which only depends on $K^*$ and $c$.
\end{teorema}

\begin{proof}
First, let us consider the integer numbers $k$ with $l_{k}+l_{k+1}
\ge r_k$. The inequality $r_k-l_k-l_{k+1}\le 2c - 2 \min\{l_k,
l_{k+1}\}$ (which is equivalent to $r_k\le 2c + |l_k-l_{k+1}|$)
gives
$$
e^{-\frac12 (l_{k}+l_{k+1}-r_k)_+} + (r_k-l_{k}-l_{k+1})_+
= e^{\frac12 (r_k-l_k-l_{k+1})}
\le e^{c - \min\{l_k, l_{k+1}\}}
\le e^{c} \big(e^{-l_k} + e^{-l_{k+1}}\big)
 \,.
$$
And now, consider the integer numbers $k$ with $l_{k}+l_{k+1} \le
r_k$. The inequality $0 \le r_k-l_k-l_{k+1}\le 2c - 2 \min\{l_k,
l_{k+1}\}$ gives $\min\{l_k, l_{k+1}\} \le c$, and consequently
$$
e^{-c}  \le e^{-\min\{l_k, l_{k+1}\}} ,
\qquad
1  \le e^{c} \big(e^{-l_k} + e^{-l_{k+1}}\big) \,.
$$
Hence
$$
e^{-\frac12 (l_{k}+l_{k+1}-r_k)_+} + (r_k-l_{k}-l_{k+1})_+
= 1 + r_k-l_{k}-l_{k+1}
\le 1 + 2c
\le (1 + 2c) \, e^{c} \big(e^{-l_k} + e^{-l_{k+1}}\big)
 \,.
$$
Then
$$
\aligned
e^{-\frac12 (l_{k}+l_{k+1}-r_k)_+} + (r_k-l_{k}-l_{k+1})_+
& \le (1 + 2c) \, e^{c} \big(e^{-l_k} + e^{-l_{k+1}}\big)  \,,
\\
e^{-l_k} + e^{-l_{k+1}} \le \Delta(k)
& \le \big( 1 + (1 + 2c) \, e^{c}\big) \big(e^{-l_k} + e^{-l_{k+1}}\big)  \,,
\endaligned
$$
for every $k\ge 1$.
Hence, if we apply Theorem \ref{t:caractren} we obtain the
conclusion, with $\inf_{m\in [A_n(h),B_n(h)]}\Gamma_{nm}^*(h)$
instead of $\min_{m\in [A_n(h),B_n(h)]}\Gamma_{nm}^*(h)$.
In order to see that the infimum is attained
we can follow an argument similar to the one at the end of the
proof of Theorem \ref{t:caractren}.
\end{proof}

\begin{teorema}
\label{t:sel}
Let us consider a train $\O$ such that there exists a constant $c>0$ with
$r_n\le c$ for every $n \ge 1$.
Then $\O$ is hyperbolic if and only if
$$
K^0:=\sup_{n\ge 1} \sup_{h\in [0,l_n]} \min_{m\in [A_n(h),B_n(h)]}
\Gamma_{nm}^0(h) < \infty\,.
$$
Furthermore,
if $\O$ is $\d$-hyperbolic, then $K^0$ is bounded by a constant
which only depends on $\d$ and $c$;
if $K^0<\infty$, then $\O$ is $\d$-hyperbolic, with $\d$ a constant
which only depends on $K^0$ and $c$.
\end{teorema}

\noindent
{\bf Remark.}
Notice that $\Gamma_{nm}^0$ is much simpler than $\Gamma_{nm}$:

Firstly, the four terms in the definition of $\D(k)$
are replaced by its first term.

Furthermore, in the first and fifth cases in the definition of $\Gamma_{nm}^0$
we remove the first term in the corresponding definition of
$\Gamma_{nm}$.

In order to obtain these simplifications, we must pay with
the hypothesis $r_n \le c$, but
this is a usual hypothesis: for instance, every flute surface satisfies it.

\begin{proof}
Notice that
$\big(r_{m} + h - l_{m+1} \big)_+ \le r_{m} \le c$
if $m=A_n(h)$ (since $l_{m+1}> h$) and
$\big(r_{m-1} + h - l_{m-1} \big)_+ \le r_{m-1} \le c$
if $m=B_n(h)$.

Hence, if we apply Theorem \ref{t:sel0} we obtain the conclusion,
with $\inf_{m\in [A_n(h),B_n(h)]}\Gamma_{nm}^0(h)$ instead of
$\min_{m\in [A_n(h),B_n(h)]}\Gamma_{nm}^0(h)$.

In order to see that the infimum is attained
we can follow an argument similar to the one at the end of the
proof of Theorem \ref{t:caractren}.
\end{proof}

\begin{prop}
\label{p:minimo2}
In any train $\O$ we have
$$
\min_{m\in [A_n(h),B_n(h)]} \Gamma_{nm}^0(h)
= \min_{m\ge 1} \Gamma_{nm}^0(h) \,,
$$
for every $n\ge 1$ and $0\le h \le l_n$.
\end{prop}

\begin{proof}
Fix $n\ge 1$ and $0\le h \le l_n$.
If $m< A_n(h)$, then
$\Gamma_{nm}^0(h) > \Gamma_{nA_n(h)}^0(h)$:
$$
\Gamma_{nm}^0(h)
\ge e^h \!\!\! \sum_{k=m+1}^{n} \!\!\! e^{-l_k}
> e^h \!\!\!\!\!\!\! \sum_{k=A_n(h)+1}^{n} \!\!\!\!\!\!\! e^{-l_k}
= \Gamma_{nA_n(h)}^0(h) \,.
$$
The case $m > B_n(h)$ is similar.
\end{proof}

Theorem \ref{t:sel} let us obtain an alternative proof of a result that appears in \cite{APRT}, but using now a completely new argument.
It is a simple sufficient condition for the hyperbolicity.

\begin{corolario}
\label{t:suficiente}
Let us consider a train $\O$ with $l_1\le l^0$, $r_n\le c_1$ for every $n$ and
\begin{equation}
\label{eq:sel}
\sum_{k=n}^{\infty} e^{-l_k} \le c_2 \, e^{-l_n} ,
\qquad \text{ for every }\, n>1\, .
\end{equation}
Then $\O$ is $\d$-hyperbolic, where $\d$ is a constant which
only depends on $c_1$, $c_2$ and $l^0$.
\end{corolario}

\noindent
{\bf Examples.}
Let us consider an increasing $C^1$ function $f$ with
$\lim_{x\to \infty} f(x)=\infty$, and define $l_n:=f(n)$ for every $n$.
A direct computation gives that $\{l_n\}$
satisfies $(\ref{eq:sel})$ if and only if there exist
constants $c, M$ with $f'(x) \ge c >0$ for every $x\ge M$.

Consequently, for $a,b>0$ and $c\in \RR$, the sequence
$l_n:=a n^b+c$ satisfies $(\ref{eq:sel})$
if and only if $b \ge 1$.

\begin{proof}
Let us consider $n\ge 1$ and $h \in [l^0,l_n]$.
Since $l_1 \le l^0 \le h$, we have that $m=A_n(h)$
satisfies $l_m \le h < l_{m+1}$ and
$$
\Gamma_{nm}^0(h)
= e^{h} \!\!\! \sum_{k=m+1}^{n} \!\!\! e^{-l_{k}}
\le e^{h} \, c_2 \, e^{-l_{m+1}}
< c_2 \, .
$$
If $h \in [0,l^0]$, then $\Gamma_{nn}^0(h) \le h \le l^0$.
Hence, $K^0 \le \max \{c_2, l^0\}$, and Theorem \ref{t:sel} gives the result.
\end{proof}

\begin{lema}
\label{l:crec0}
$\,$

\begin{enumerate}
\item[$(1)$]
Let us consider a sequence $\{l_n\}$ such that $l_{m}\le l_n+c$
for every positive integer numbers $m\le n$.
Then there exists a non-decreasing sequence $\{l'_n\}$,
such that $|l_n-l'_n| \le c$ for every $n$.

\item[$(2)$]
Let us consider a non-decreasing sequence $\{l'_n\}$.
If $\{l_n\}$ is a sequence with $|l_n-l'_n| \le c$ for every $n$,
then $l_{m}\le l_n+2\, c$ for every positive integer numbers $m\le n$.
\end{enumerate}
\end{lema}

\begin{proof}
We prove now the first part of the lemma.
We define a sequence $\{l'_n\}$ in the following way:
$l'_{n}:=\max\{l_1,l_2,\dots,l_n\}$.
It is clear that $\{l'_n\}$ is a non-decreasing sequence.
Since $l_{m}\le l_n+c$ for every $m=1,2,\dots,n,$ we have
$l_n \le l'_n \le l_n+c$.
Consequently, $|l_n-l'_n| \le c$ for every $n$.

In order to prove the second part, notice that if
$m\le n$, then $l_{m} \le l'_m+c \le l'_n+c\le l_n+2\,c$.
\end{proof}

The two following theorems provide necessary conditions for hyperbolicity.

\begin{teorema}
\label{t:necesario1}
Let us consider an hyperbolic train $\O$ with
$l_{m}\le l_n+c_1$ for every positive integer numbers $m\le n$.
If $K$ is the constant defined in Theorem \ref{t:caractren},
then
$$
r_n \le 2 \max\{K,1\} + 2 \log \max\{K,1\} + 3 \,c_1 \,,
\qquad \text{ for every $n$ with }\,
l_{n+1}> 4(K+c_1)\,.
$$
\end{teorema}

\begin{proof}
Let us define $M:= \max\{K,1\}$ and fix $n$ with $l_{n+1}> 4(K+c_1)$.

Let us assume that $r_n \le l_{n+1}$. Consider $\varepsilon \in (0, 1/2)$
and $h_{n+1}:= l_{n+1}- \varepsilon r_n$. Then
$$
\aligned
\Gamma_{n+1,n+1}(h_{n+1}) & = \min\{ l_{n+1}- \varepsilon r_n,\, \varepsilon r_n\}
= \varepsilon r_n\,,
\\
\Gamma_{n+1,m}(h_{n+1}) & \ge l_{m}- h_{n+1}
\ge l_{n+1}- c_1 - h_{n+1}
=\varepsilon r_n- c_1\,,
\qquad \quad \text{if $m>n+1$}\,,
\\
\Gamma_{n+1,n}(h_{n+1}) & \ge (r_{n}+ h_{n+1} -l_{n+1})_+
=(1-\varepsilon) r_n \,,
\qquad \quad \text{if $l_n \le h_{n+1}$}\,,
\\
\Gamma_{n+1,m}(h_{n+1}) & \ge e^{h_{n+1}} \Delta(n)
\ge e^{l_{n+1}- \varepsilon r_n} e^{-\frac12 (l_{n} + l_{n+1}- r_n)}
\ge e^{l_{n+1}- \varepsilon r_n} e^{-\frac12 (l_{n+1} + l_{n+1} + c_1 - r_n)}
\\
& = e^{-\frac12 c_1 + (\frac12 - \varepsilon) r_n} ,
\qquad \quad \text{if either $m<n$ or $m=n$ and $l_n > h_{n+1}$}\,.
\endaligned
$$
Since $\varepsilon \in (0, 1/2)$
$$
M \ge \min \big\{ \varepsilon r_n,\, \varepsilon r_n- c_1,\,
(1-\varepsilon) r_n,\, e^{-\frac12 c_1 + (\frac12 - \varepsilon) r_n}\big\}
= \min \big\{ \varepsilon r_n- c_1,\, e^{-\frac12 c_1 + (\frac12 - \varepsilon) r_n}\big\}\,,
$$
and we deduce
$$
r_n \le \max \Big\{ \frac{M+c_1}{\varepsilon}\,, \,
\frac{\log M+c_1/2}{1/2 - \varepsilon}\Big\}\,.
$$
Taking $\varepsilon = (M+c_1)/(2M+2\log M+3c_1)$
(notice that $\varepsilon \in (0,1/2)$, since $\log M \ge 0$), we obtain
the equality of the two terms inside the maximum, and therefore
$r_n \le 2M+2\log M+3c_1$.

We prove now that $r_n \le l_{n+1}$.
Seeking for a contradiction, assume that
$r_n > l_{n+1}$, and consider $h^{n+1}:= \frac34 \,l_{n+1}$.
A similar argument, with $h^{n+1}$ instead of $h_{n+1}$, gives:

\noindent
If $l_{n} + l_{n+1} < r_n$, since $l_{n+1}> 4(K+c_1)$,
$$
K \ge \min \Big\{\frac14\, l_{n+1} ,\, \frac14\, l_{n+1} - c_1 ,\,
\frac34\, l_{n+1},\, e^{\frac34 l_{n+1}}\Big\}
= \frac14\, l_{n+1} - c_1 > K ,
$$
since $l_{n+1}> 4(K+c_1)$, and this is a contradiction.
If $l_{n} + l_{n+1} \ge r_n$, we obtain with a similar argument
$$
K \ge \min \Big\{\frac14\, l_{n+1} ,\, \frac14\, l_{n+1} - c_1 ,\,
\frac34\, l_{n+1},\, e^{\frac14 l_{n+1} -\frac12 c_1}\Big\}
= \min \Big\{ \frac14\, l_{n+1} - c_1 ,\,
 e^{\frac14 l_{n+1} -\frac12 c_1}\Big\} > K ,
$$
since $l_{n+1}> 4(K+c_1)$, and this is the contradiction we are looking for.
\end{proof}

Condition $l_{m}\le l_n+c_1$ for every positive integer numbers
$m\le n$ in Theorem \ref{t:necesario1} can seem superfluous, but we
have examples which prove that, in fact, if it is removed, then the
conclusion of the theorem is not true.

The following theorem obtains a similar inequality to (\ref{eq:sel}) but with an explicit control of the constants involved.

\begin{teorema}
\label{t:necesario2}
Let us consider an hyperbolic train $\O$ with
$l_{m}\le l_n+c_1$ for every positive integer numbers $m\le n$.
If $K$ is the constant defined in Theorem \ref{t:caractren},
then
$$
\sum_{k=n}^{\infty} e^{-l_k} \le K \, e^{K+c_1} \, e^{-l_n}  ,
\qquad \text{ for every $n$ with }\,
l_n > 2 K + c_1 \, .
$$
\end{teorema}

\begin{proof}
Theorem \ref{t:caractren}
and Proposition \ref{p:minimo}
give that
$$
\min_{m \ge 1}\Gamma_{n m}(h) \le K\,,
\qquad \text{ for every $n \ge 1$ and }\,
h \in [0, l_n] \,.
$$
Let us fix $n$ with $l_n > 2 K + c_1$ and $n_0 \ge n$.
Consider $\varepsilon > 0$ with
$l_n \ge 2 K + c_1 + \varepsilon$.
If we define $h:= l_n - K - c_1 - \varepsilon/2 \ge K + \varepsilon/2 > K$,
then for any $m \ge n$ we have
$l_m - h \ge l_n - h - c_1 = K + \varepsilon/2 > K$ and
$$
\Gamma_{n_0 m}(h) \ge \Gamma_{n_0 m}^0 (h) \ge K + \varepsilon/2 > K \,.
$$
If $m < n$, we obtain
$$
\Gamma_{n_0 m}(h) \ge \Gamma_{n_0 m}^0 (h) \ge e^h \sum_{k=n}^{n_0} e^{-l_k} .
$$
Consequently,
$$
K \ge \min_{m \ge 1}\Gamma_{n_0 m}(h)
= \min_{1 \le m < n}\Gamma_{n_0 m}(h)
\ge e^{l_n - K - c_1 - \varepsilon/2} \sum_{k=n}^{n_0} e^{-l_k} ,
$$
for every $n_0 \ge n$ and $\varepsilon$ small enough. Therefore
$$
K \ge e^{l_n - K - c_1} \sum_{k=n}^{\infty} e^{-l_k} ,
$$
which finishes the proof.
\end{proof}

The last three theorems, Theorem \ref{t:caractren} and Proposition
\ref{p:lacotado} give the following powerful and simple characterization.
In particular, this result characterizes hyperbolicity of trains for which $l_n$ is a non-decreasing sequence.

\begin{teorema}
\label{t:caracsel}
Let us consider a train $\O$ with
$l_{m}\le l_n+c_1$ for every positive integer numbers $m\le n$.

\begin{enumerate}
\item[$(1)$]
If $\{l_n\}$ is a bounded sequence,
then $\O$ is hyperbolic.

\item[$(2)$]
If $\lim_{n\to \infty} l_n = \infty$,
then $\O$ is hyperbolic if and only if
$\{r_n\}$ is a bounded sequence and
$(\ref{eq:sel})$ holds for some constant $c_2$.
\end{enumerate}
\end{teorema}

\noindent
{\bf Remark.}
Note that Theorem \ref{t:caracsel}
deals with every case under the hypothesis
``$l_{m}\le l_n+c_1$ for $m\le n$":
$\{l_n\}$ is either a bounded sequence or a sequence with limit $\infty$.

\medskip

If we have an hyperbolic train, we want to study
what kind of transformations in $\{l_n\}$ and $\{r_n\}$
allows to obtain another hyperbolic train.

\begin{teorema}
\label{t:sucesiones}
Consider two trains $\O$ and $\O'$.
Let us assume that $\O$ is $\d$-hyperbolic.
Then, $\O'$ is $\d'$-hyperbolic if we have either:

\begin{enumerate}
\item[$(1)$]
$l'_n = l_n$ and
$r'_n \le r_n$ for every $n$
$($and then $K' \le K)$, or

\item[$(2)$]
$l'_n =\l l_n$
and
$r'_n =\l r_n$ for every $n$
$(\l \ge 1)$
$($and then $K' \le \l K + (1+\l) K^{\l})$, or

\item[$(3)$]
$l'_n =\l l_n$
and
$r'_n =\mu r_n$ for every $n$
$(\l \ge 1,\, 0 \le \mu \le \l)$
$($and then $K' \le \l K + (1+\l) K^{\l})$.
\end{enumerate}
\end{teorema}

\begin{proof}
In case $(1)$, $\big(\Gamma_{nm}\big)'(h) \le \Gamma_{nm}(h)$ for every $n,m\ge 1$, since
$\Gamma_{nm}(h)$ is a non-decreasing function in each variable $r_k$.
This allows to deduce $(1)$.

In order to prove the second part, notice that
(since $\l \ge 1$)
$$
e^{\l h} \displaystyle\sum_{k} \big( e^{-\l l_{k}} + e^{-\l l_{k+1}} +
e^{-\frac12 (\l l_{k}+\l l_{k+1}-\l r_k)_+} \big)
\le
\Big( e^{h} \displaystyle\sum_{k} \big( e^{-l_{k}} + e^{-l_{k+1}} +
e^{-\frac12 (l_{k}+l_{k+1}-r_k)_+} \big)
\Big)^{\l}
\,.
$$
Notice that $t \le (1+t)^{\lambda}$
for every $t\ge 0$ and $\l \ge 1$. Hence, if $r_k-l_{k}-l_{k+1}\ge 0$,
$$
\aligned
e^{\l h} \displaystyle\sum_{k} (\l r_k-\l l_{k}-\l l_{k+1})_+
& \le
\l \, e^{\l h} \displaystyle\sum_{k} \big( 1 + (r_k-l_{k}-l_{k+1})_+ \big)^{\l}
\\
& \le
\l \Big( e^{h} \displaystyle\sum_{k} \big( e^{-\frac12 (l_{k}+l_{k+1}-r_k)_+} + (r_k-l_{k}-l_{k+1})_+ \big)
\Big)^{\l}
\,.
\endaligned
$$
We also have
$$
\big(\l r_{m} + \l h - \l l_{m+1} \big)_+ = \l \big(r_{m} + h - l_{m+1} \big)_+ \,,
\quad
\l l_m - \l h = \l (l_m - h) \,,
\quad
\min  \big\{\l h,\, \l l_n - \l h \big\} = \l \min  \big\{h,\, l_n - h \big\} \,.
$$
Consequently,
$\big(\Gamma_{nm}\big)'(\l h) \le \l \Gamma_{nm}(h) + \Gamma_{nm}(h)^{\l}
+ \l \Gamma_{nm}(h)^{\l}$ for every $n,m\ge 1$ and $0 \le h \le l_n$,
and then $K' \le \l K + (1+\l) K^{\l}$.

Item $(3)$ is a direct consequence of $(1)$ and $(2)$.
\end{proof}

We want to study now the following question:
If we have an hyperbolic train with $\{r_n\}\in l^{\infty}$,
what kind of perturbations are allowed on $\{l_n\}$
so that the train is still hyperbolic?
Theorem \ref{t:union} answers this question providing a
great deal of hyperbolic flute surfaces.

We need the following definitions.

\begin{definicion}
We denote by $H$ the following set of sequences:
$$
\aligned
H
:\!&=\,\big\{ \{x_n\}: \, \text{ the train with }
l_n=x_n \, \text{ and }
r_n=0 \, \text{ for every $n$ is hyperbolic} \, \big\}
\\
&=\, \big\{ \{x_n\}: \, \text{ every train with }
l_n=x_n \, \text{ for every $n$ and }
\{r_n\} \in l^{\infty} \, \text{ is hyperbolic} \, \big\} \,.
\endaligned
$$
\end{definicion}

The second equality is a direct consequence of Theorem
\ref{t:HypStab}.

\begin{definicion}
We say that the sequence $\{y_n\}$ is a union of
the sequences $\{x^1_n\}, \dots ,\{x^N_n\},$ if
$\{x^1_n\}, \dots ,\{x^N_n\}$ are subsequences of $\{y_n\}$, and
$\{x^1_n\}, \dots ,\{x^N_n\}$ is a partition of $\{y_n\}$.
\end{definicion}

\begin{teorema}
\label{t:union}
Let us consider a sequence $\{l_n\} \in H$.

\begin{enumerate}
\item[$(1)$]
If $l'_n=l_n + x_n$ with $\{x_n\} \in l^{\infty}$,
then $\{l'_n\} \in H$.

\item[$(2)$]
Fix a positive integer $N$.
Let us assume that $\{l_n\}$ is a subsequence $\{l'_{n_k}\}$ of $\{l'_n\}$
such that $n_{k+1} - n_k \le N$ for every $k$,
and $\max \{l'_{n_{k}}, l'_{n_{k+1}} \} \le l'_m +N$
for every $m\in (n_{k}, n_{k+1})$ and every $k$.
Then $\{l'_n\} \in H$.

\item[$(3)$]
If $\{l'_n\}$ is any union of
the sequences $\{l^1_n\}, \dots ,\{l^N_n\} \in H$,
then $\{l'_n\} \in H$.

\item[$(4)$]
If $\{l'_n\}$ is a union of $\{l_n\}$ and
a sequence $\{x_n\} \in l^{\infty}$,
then $\{l'_n\} \in H$.

\item[$(5)$]
Let us assume that $\{l'_n\}$ is any union of
the sequences $\{l^1_n\}, \dots ,\{l^N_n\}$
which verify
$$
\sum_{k=n}^{\infty} e^{-l^j_k} \le c \, e^{-l^j_n}  ,
\qquad \text{ for every $n>1$ and } \, j=1,\dots,N\, .
$$
Then $\{l'_n\} \in H$.

\item[$(6)$]
Fix a positive integer $N$.
Let us assume that $\{x_n\}$ is a subsequence $\{l'_{n_k}\}$ of $\{l'_n\}$
such that $\max \{l'_{n_{k}}, l'_{n_{k+1}} \} \le l'_m +N$
for every $m\in (n_{k}, n_{k+1})$ and every $k$.
If $\{x_n\} \notin H$, then $\{l'_n\} \notin H$.

\item[$(7)$]
Fix a positive integer $N$. Let $\s$ be a permutation of the
positive integer numbers such that $|\s(n)-n| \le N$ for every $n$,
and consider $l'_n:=l_{\s(n)}$. Then $\{l'_n\} \in H$.
\end{enumerate}
\end{teorema}

\noindent
{\bf Remarks.}

$(1)$ In fact, $(7)$ gives the following stronger statement: If $\s$
is a permutation of the positive integer numbers such that
$|\s(n)-n| \le N$ for every $n$, then $\{l_{\s(n)}\} \in H$ if and
only if $\{l_n\} \in H$ (since $\s^{-1}$ also satisfies
$|\s^{-1}(n)-n| \le N$ for every $n$).

$(2)$ We have examples showing that the conclusions of Theorem
\ref{t:union} do not hold if we remove any of the hypothesis.

\begin{proof}
$(1)$ is a direct consequence of Theorem \ref{t:HypStab}.

\spb

$(2)$
Fix $n\ge 1$ and $h\in [0, l'_n]$.

Let us consider the maximum integer $k_0$
such that $n_{k_0} \le n < n_{k_0+1}$.

If $l'_{s} \le h$ for some $s \in [n_{k_0}, n_{k_0+1}]$,
by symmetry, without loss of generality we can assume that
there exists some $s\in [n_{k_0}, n)$ with $l'_{s} \le h$
(the case $s=n$ is trivial: if $l'_n \le h$, then
$h=l'_n$ and $\big( \Gamma_{nn}^0 \big)'(h)=0$).
Hence $A'_n(h)\in [n_{k_0}, n)$ and
then $l'_{k} \ge h$ for every $k\in (A'_n(h),n]$ and
$n-A'_n(h) \le n-n_{k_0} \le N-1$; consequently,
$$
\big( \Gamma_{n A'_n(h)}^0 \big)'(h)
= \!\!\!\!\! \sum_{k=A'_n(h)+1}^{n}  \!\!\!\!\!\!\! e^{h-l'_{k}}
\le \!\!\!\! \sum_{k=A'_n(h)+1}^{n} \!\!\!\!\!\! 1
= n - A'_n(h)
\le N-1 \,.
$$
\indent
Let us assume now that $l'_{s} > h$ for every $s \in [n_{k_0}, n_{k_0+1}]$.
There exists some integer $m$ with $\Gamma_{k_0 m}^0 (h) \le K^0$.
By symmetry, without loss of generality we can assume that $m \le k_0$.

If $m=k_0$, then $\min \{ h, \, l_{k_0}-h \} \le K^0$.
If $\min \{ h, \, l_{k_0}-h \}= h$,
then $h \le K^0$ and we can deduce
$$
\big( \Gamma_{nn}^0 \big)'(h) = \min \{ h, \, l'_{n}-h \}
\le h \le K^0 .
$$
If $\min \{ h, \, l_{k_0}-h \}= l_{k_0}-h$, then $l_{k_0}-h \le K^0$ and
$$
\big( \Gamma_{n n_{k_0}}^0 \big)'(h)
= l'_{n_{k_0}}-h + \sum_{k=n_{k_0}}^{n} e^{h-l'_{k}}
\le l_{k_0}-h + \sum_{k=n_{k_0}}^{n} 1
\le K^0 + N .
$$
\indent
If $m< k_0$ and $l_m>h$, then $\Gamma_{{k_0}m}^0 (h)
= l_m -h + e^h \sum_{k=m}^{k_0} e^{-l_{k}} \le K^0$.
Hence
$$
\aligned
\big( \Gamma_{n n_{m}}^0 \big)'(h)
& = l'_{n_{m}}-h + e^h \sum_{k=n_{m}}^{n_{k_0}} e^{-l'_{k}}
+ \sum_{k=n_{k_0}+1}^{n} e^{h-l'_{k}}
\\
& \le l'_{n_{m}}-h + e^h \Big( e^{-l'_{n_m}} +
\sum_{j=m+1}^{k_0} \sum_{k=n_{j-1}+1}^{n_{j}} e^{-l'_{k}} \Big)
+ \sum_{k=n_{k_0}+1}^{n} 1
\\
& \le l'_{n_{m}}-h + e^h \Big( e^{-l'_{n_m}} +
\sum_{j=m+1}^{k_0} N \, e^{N-l'_{n_j}} \Big) + N - 1
\\
& \le N \, e^N \Big( l_{m}-h + e^h \sum_{j=m}^{k_0} e^{-l_{j}} \Big) + N - 1
\le N \, e^N \, K^0 + N - 1 \,.
\endaligned
$$

If $m< k_0$ and $l_m\le h$, a similar argument gives the same bound for
$\big( \Gamma_{n n_{m}}^0 \big)'(h)$.

Then, $\big( K^0 \big)' \le N \, e^N \, K^0 + N$
and Theorem \ref{t:sel} implies $(2)$.

\spb

$(3)$
Assume first that $N=2$; then
$\{l'_n\}$ is the union of $\{l^1_n\}$ and $\{l^2_n\}$.
We denote by $\{l'_{n_k^i}\}$
the subsequence $\{l^i_n\}$ in $\{l'_n\}$, for $i=1,2$.
Fix $n\ge 1$ and $h\in [0, l'_n]$.
By symmetry,
without loss of generality we can assume that
there exist $k_1$ with $n_{k_1}^1=n$ and $m_1 \le k_1$
with $\big(\Gamma_{k_1 m_1}^0 \big)^1 (h) \le (K^{0})^1$.

We can assume that $l'_{s} > h$ for every
$s \in (n_{m_1}^1, n_{k_1}^1)$,
since the other case is similar.

If there is no $k$ with $n_k^2 \in [n_{m_1}^1, n_{k_1}^1]$,
then
$\big(\Gamma_{n_{k_1}^1 n_{m_1}^1}^0 \big)' (h)
= \big(\Gamma_{k_1 m_1}^0 \big)^1 (h) \le (K^{0})^1$.

Assume now that there exists $k$ with $n_k^2 \in (n_{m_1}^1, n_{k_1}^1)$.
Let us define
$k_2 := \max \{k: \, n_k^2 \in (n_{m_1}^1, n_{k_1}^1)\}$.

If there exists $m_2 \le k_2$ such that
$\big(\Gamma_{k_2 m_2}^0 \big)^2 (h) \le (K^{0})^2$, then
$$
\big(\Gamma_{n_{k_1}^1, \, \max\{n_{m_1}^1, \, n_{m_2}^2\}}^0 \big)' (h)
\le \big(\Gamma_{k_1 m_1}^0 \big)^1 (h)
+ \big(\Gamma_{k_2 m_2}^0 \big)^2 (h) \le (K^{0})^1 + (K^{0})^2 .
$$
\indent If there exists $k_3$ verifying the next three conditions
simultaneously:

$(a)$ $n_{k_3}^2 \in (n_{m_1}^1, n_{k_1}^1)$,

$(b)$ there exists $m_3 \le k_3$ such that
$\big(\Gamma_{k_3 m_3}^0 \big)^2 (h) \le (K^{0})^2$,

$(c)$ for every $k \in (k_3,k_2]$ we have
$\big(\Gamma_{k m}^0 \big)^2 (h) > (K^{0})^2$ for every $m \le k$,

\noindent
then there exists $m_0 > k_2$ such that
$\big(\Gamma_{k_3+1, m_0}^0 \big)^2 (h) \le (K^{0})^2$:
In fact, seeking for a contradiction, let us assume that
there exists $m_0 \in (k_3+1,k_2]$ with
$\big(\Gamma_{k_3+1, m_0}^0 \big)^2 (h) \le (K^{0})^2$;
then
$\big(\Gamma_{m_0 m_0}^0 \big)^2 (h)
\le \big(\Gamma_{k_3+1, m_0}^0 \big)^2 (h) \le (K^{0})^2$
(recall that $l'_{s} > h$ for every $s \in (n_{m_1}^1, n_{k_1}^1)$),
which is actually a contradiction with $(c)$.
Hence,
$$
\big(\Gamma_{n_{k_1}^1, \, \max\{n_{m_1}^1, \, n_{m_3}^2\}}^0 \big)' (h)
\le \big(\Gamma_{k_1 m_1}^0 \big)^1 (h)
+ \big(\Gamma_{k_3 m_3}^0 \big)^2 (h)
+ \big(\Gamma_{k_3+1, m_0}^0 \big)^2 (h) \le (K^{0})^1 + 2 (K^{0})^2 .
$$
\indent
If for any $k$ with $n_{k}^2 \in (n_{m_1}^1 n_{k_1}^1)$ we have
$\big(\Gamma_{k m}^0 \big)^2 (h) > (K^{0})^2$ for every $m \le k$,
let us define
$k_4 := \min \{k: \, n_k^2 \in (n_{m_1}^1, n_{k_1}^1)\}$.
As in the last case,
then there exists $m_4 > k_2$ such that
$\big(\Gamma_{k_4 m_4}^0 \big)^2 (h) \le (K^{0})^2$, and hence
$$
\big(\Gamma_{n_{k_1}^1 n_{m_1}^1}^0 \big)' (h)
\le \big(\Gamma_{k_1 m_1}^0 \big)^1 (h)
+ \big(\Gamma_{k_4 m_4}^0 \big)^2 (h)
\le (K^{0})^1 + (K^{0})^2 .
$$
\indent
Consequently,
$\big(K^0 \big)' \le 2 (K^{0})^1 + 2 (K^{0})^2$
and Theorem \ref{t:sel} implies $(3)$ with $N=2$.
The result for $N$ sequences is obtained by applying
$N-1$ times this result for $2$ sequences.

\spb

$(4)$ is a direct consequence of $(3)$ and Proposition \ref{p:lacotado}.

\spb

$(5)$ is a direct consequence of $(3)$ and Theorem \ref{t:suficiente}.

\spb

$(6)$
Since $\{x_n\} \notin H$,
by Theorem \ref{t:sel}
and Proposition \ref{p:minimo2},
for each $M > N$ there exist $k_0$ and $h\in (0, x_{k_0})$ with
$\Gamma_{{k_0}m}^0(h) \ge M$,
for every $m \ge 1$.

Consider $m \ge 1$.
By symmetry, without loss of generality we can assume that
$m \le n_{k_0}$. If $m = n_{k_0}$, then
$$
\big(\Gamma_{n_{k_0}n_{k_0}}^0\big)' (h)
= \min  \big\{h,\, l'_{n_{k_0}} - h \big\}
= \min  \big\{h,\, x_{k_0} - h \big\}
= \Gamma_{{k_0}k_0}^0(h) \ge M .
$$
Notice that if $m \in (n_{k_0-1}, n_{k_0})$, then
$$
l'_m - h
\ge l'_{n_{k_0}} - h - N
=  x_{k_0} - h - N
\ge \Gamma_{k_0k_0}^0 (h) - N
\ge M - N > 0 \,,
$$
and $l'_m > h$. Hence
$\big(\Gamma_{n_{k_0} m}^0\big)' (h) \ge l'_m - h \ge M - N$.

In the case $m \le n_{k_0-1}$, we have
$n_{k_1-1} < m \le n_{k_1}$ for some $k_1 < k_0$.

If $x_{k_1} \le h$, then
$$
\big(\Gamma_{n_{k_0} m}^0\big)' (h)
\ge e^{h} \! \displaystyle\sum_{k=m+1}^{n_{k_0}} \! e^{-l'_{k}}
\ge e^{h} \! \displaystyle\sum_{k=k_1+1}^{k_0} \! e^{-x_{k}}
= \Gamma_{k_0k_1}^0(h)
\ge M .
$$
\indent
If $x_{k_1} > h$ and $l'_m>h$, then
$$
\aligned
\big(\Gamma_{n_{k_0} m}^0\big)' (h)
& = l'_m - h +
e^{h} \displaystyle\sum_{k=m}^{n_{k_0}} e^{-l'_{k}}
\ge l'_{n_{k_1}} - h - N +
e^{h} \displaystyle\sum_{k=m}^{n_{k_0}} e^{-l'_{k}}
\\
& \ge  x_{k_1} - h - N +
e^{h} \displaystyle\sum_{k=k_1}^{k_0} e^{-x_{k}}
= \Gamma_{k_0k_1}^0 (h) - N
\ge M - N .
\endaligned
$$
\indent
If $x_{k_1} > h$ and $l'_m \le h$, then
$x_{k_1} - N = l'_{n_{k_1}} - N \le l'_m \le h$ and
$0 \ge x_{k_1} - h - N$; therefore
$$
\big(\Gamma_{n_{k_0} m}^0\big)' (h)
= e^{h} \! \displaystyle\sum_{k=m+1}^{n_{k_0}} \! e^{-l'_{k}}
\ge x_{k_1} - h - N
+ e^{h} e^{-x_{k_1}} -1
+ e^{h} \! \displaystyle\sum_{k=k_1+1}^{k_0} \! e^{-x_{k}}
= \Gamma_{k_0k_1}^0(h) - N - 1
\ge M - N - 1 .
$$
Consequently,
$\big(K^0\big)' \ge M - N - 1$ for every $M > N$,
and hence $\big(K^0\big)' = \infty$.
Then $\{l'_n\} \notin H$ by Theorem \ref{t:sel}.

\spb

$(7)$
First, we want to remark the following elementary fact:
If $i<j$ and $\s(i)>\s(j)$, then $|i-j|< 2N$:
$|i-j|=j-i<j-\s(j) +\s(i)-i \le 2N$.

Fix $n\ge 1$ and $h\in [0, l'_n]$.
There exists $\s(m)$ with $\Gamma_{\s(n)\s(m)}^0(h) \le K^0$.
By symmetry, without loss of generality we can assume that
$\s(m)\le \s(n)$.

If $m=n$, then $\s(m)=\s(n)$ and
$\big(\Gamma_{nn}^0\big)'(h) =
\Gamma_{\s(n)\s(n)}^0(h) \le K^0$.

We consider now the case
$\s(m)< \s(n)$.

If $m>n$, then $m-n < 2N$.

$\quad$
If $B'_n(h) > m$, then $l'_k > h$ for every $k \in (n,m]$ and
$$
\big(\Gamma_{nm}^0\big)'(h)
= l'_m - h + \sum_{k=n}^{m} e^{h-l'_{k}}
\le l_{\s(m)} - h + 2N
\le \Gamma_{\s(n)\s(m)}^0(h) + 2N
\le K^0 + 2N .
$$
\indent
$\quad$
If $B'_n(h) \le m$, then $l'_k > h$ for every $k \in (n,B'_n(h))$ and
$$
\big(\Gamma_{n B'_n(h)}^0\big)'(h)
= \!\! \sum_{k=n}^{B'_n(h)-1} \!\!\! e^{h-l'_{k}}
\le 2N .
$$
\indent
We deal now with the case $m<n$.
Notice first that $\s([m,n])\subset [m-N,n+N]$
and $[m+N,n-N]\subset [\s(m),\s(n)]$;
then, in $\s([m,n])\setminus [\s(m),\s(n)]$
there are at most $4N$ integers.

If $A'_n(h) \ge m$, then $l'_k > h$ for every $k \in (A'_n(h),n)$, and
$$
\aligned
\big(\Gamma_{n A'_n(h)}^0\big)'(h)
& = e^{h} \!\!\!\!\!\! \sum_{k=A'_n(h)+1}^n \!\!\!\!\!\! e^{-l'_{k}}
\le e^{h} \!\!\!  \sum_{\substack{k \in [m,n] \\ l_{\s(k)} \ge h}} \!\!\!  e^{-l_{\s(k)}}
= e^{h} \!\!\!\!\!\!  \sum_{\substack{j \in \s([m,n]) \\ l_j \ge h}} \!\!\!\!\!\!  e^{-l_j}
\le  \!\!\!\!\!\!   \sum_{\substack{j \in \s([m,n]) \setminus [\s(m),\s(n)] \\ l_j \ge h}}
\!\!\!\!\!\!\!\!\!\!\!\!\!\!  e^{h-l_j}
 + e^{h} \!\!\!  \sum_{\substack{j = \s(m) \\ l_j \ge h}}^{\s(n)} \!\!\! e^{-l_j}
\\
&
\le  4N + 1
 + e^{h} \!\!\!\!  \sum_{j = \s(m) +1}^{\s(n)} \!\!\!\! e^{-l_j}
\le  4N + 1 + \Gamma_{\s(n)\s(m)}^0(h)
\le 4N + 1 + K^{0} .
\endaligned
$$
\indent
If $A'_n(h) < m$, then $l'_k > h$ for every $k \in [m,n)$, and
$$
\aligned
\big(\Gamma_{n m}^0\big)'(h)
& = l'_m - h +  e^{h}  \sum_{k=m}^n e^{-l'_{k}}
= l_{\s(m)} - h + e^{h} \!\! \sum_{k \in [m,n]} \!\! e^{-l_{\s(k)}}
= l_{\s(m)} - h + e^{h} \!\!\!\!\! \sum_{j \in \s([m,n])} \!\!\!\!\! e^{-l_j}
\\
&
\le  \!\!\!\!\!\! \sum_{j \in \s([m,n]) \setminus [\s(m),\s(n)]}
\!\!\!\!\!\!\!\!\!\!\!\!\!\! e^{h-l_j}
 + l_{\s(m)} - h + e^{h} \!\! \sum_{j = \s(m)}^{\s(n)} \!\! e^{-l_j}
\le  4N + \Gamma_{\s(n)\s(m)}^0(h)
\le 4N + K^{0} .
\endaligned
$$
\indent
Hence, $\big(K^{0}\big)' \le 4N + 1 + K^{0}$,
and Theorem \ref{t:sel} gives $(7)$.
\end{proof}

\bigskip

\section{Trigonometric lemmas.}

\medskip

In this section some technical lemmas are collected. All of them
have been used in Section \ref{s:MainResults} in order to simplify the proof of Theorem
\ref{t:caractren}.

\begin{definicion}
Given a surface $M$, a geodesic $\g$ in $M$,
and a continuous unit vector field $\xi$ along $\g$,
orthogonal to $\g$, we define the
\emph{Fermi coordinates} based on $\g$ as the map
$E(u,v):=\exp_{\g(u)} v \xi(u)$.
\end{definicion}

It is well known that the Riemannian metric can be expressed in
Fermi coordinates as $ds^2= dv^2 + \eta^2(u,v) \, du^2$, where
$\eta(u,v)$ is the solution of the scalar equation $\p^2\eta/\p v^2
+ K\eta=0$, $\eta(u,0)=1$, $\p \eta/\p v(u,0)=0$, and $K$ is the
curvature of $M$ (see e.g. \cite[p. 247]{C}). Consequently, if
$M$ is a non-exceptional Riemann surface,
the Poincar\'e metric in Fermi coordinates (based on any
geodesic $\g$) is $ds^2= dv^2 + \cosh\!^2 v\, du^2$,
since $K=-1$ in the Poincar\'e metric.
We always consider in a train the Fermi coordinates
based on $(a_0,b_0)$.

\begin{definicion}
Let us consider Fermi coordinates $(u,v)$ in $\DD$.
We define the distances $d_1\big((u_1,v_1),(u_2,v_2)\big)$,
$d_2\big((u_1,v_1),(u_2,v_2)\big)$ as follows:
without loss of generality we can assume that $v_1\ge v_2$; then
$$
\aligned
d_1\big((u_1,v_1),(u_2,v_2)\big)
:\!&=\,d\big((u_1,v_1),(u_1,v_2)\big) + d\big((u_1,v_2),(u_2,v_2)\big)
=\,v_1-v_2 + d\big((u_1,v_2),(u_2,v_2)\big) ,
\\
d_2\big((u_1,v_1),(u_2,v_2)\big)
:\!&=\,d\big((u_1,v_1),(u_2,v_1)\big) + d\big((u_2,v_1),(u_2,v_2)\big)
=\, d\big((u_1,v_1),(u_2,v_1)\big) + v_1-v_2\,.
\endaligned
$$
\end{definicion}

The following lemma shows that the ``cartesian distances"
$d_1$ and $d_2$ are comparable to $d$.

\begin{lema}
\label{l:cartesiana}
Let us consider Fermi coordinates $(u,v)$ in $\DD$ and
the distances $d_1$ and $d_2$. Then
$$
\frac12\, d_1 \le d \le d_1\,, \qquad
\frac13\, d_2 \le d \le d_2\,.
$$
\end{lema}

\begin{proof}
Triangle inequality gives directly $d \le d_1$ and $d \le d_2$.
Let us consider $v_1\ge v_2$.
It is easy to check that
$$
d\big((u_1,v_1),(u_1,v_2)\big) \le d\big((u_1,v_1),(u_2,v_2)\big), \qquad
d\big((u_1,v_2),(u_2,v_2)\big) \le d\big((u_1,v_1),(u_2,v_2)\big)
$$
and this implies $d_1\le 2 d$.

We also have $d\big((u_2,v_1),(u_2,v_2)\big) \le d\big((u_1,v_1),(u_2,v_2)\big)$,
and then
$$
\aligned
d\big((u_1,v_1),(u_2,v_1)\big)
& \le  d\big((u_1,v_1),(u_2,v_2)\big) + d\big((u_2,v_1),(u_2,v_2)\big)
\le  2\, d\big((u_1,v_1),(u_2,v_2)\big) ,
\\
d_2\big((u_1,v_1),(u_2,v_2)\big)
& =d\big((u_1,v_1),(u_2,v_1)\big) + d\big((u_2,v_1),(u_2,v_2)\big)
\le 3 \, d\big((u_1,v_1),(u_2,v_2)\big) .
\endaligned
$$
\end{proof}

\begin{lema}
\label{l:cartesiana2}
Let $\O$ be a train and $l_0$ any positive constant.
We have
$$
d_1(z, \g_n \cap (a_n,b_n)) \le 2\, d_\O(z, (a_n,b_n))
+ 2 \Arcsinh \frac{1}{\sqrt{2 \tanh l_0}}\;,
$$
for every $n>0$ and $z\in \O$ with $l_0 \le h(z) \le l_n$.
\end{lema}

\begin{proof}
Let $w$ be the nearest point in $(a_n,b_n)$ to $z$, and define
$v:=\g_n \cap (a_n,b_n)$, let $v_0$ be the nearest point in
$(a_0,b_0)$ to $v$ and $w_0$ the nearest point in $(a_0,b_0)$ to
$w$. Consider the geodesic quadrilateral in $\O^+$ with vertices
$v$, $w$, $w_0$ and $v_0$. Standard hyperbolic trigonometry
gives that
$$
\tanh d_\O(w,w_0) = \tanh d_\O(v,v_0)  \cosh d_\O(v_0,w_0)
= \tanh l_n  \cosh d_\O(v_0,w_0)\,.
$$
Denote by $v'$ (respectively $w'$)
the point in $\g_n^+=[v,v_0]\subset \O^+$ (respectively in $[w,w_0]\subset \O^+$)
with $h(v')=h(z)$ (respectively $h(w')=h(z)$).
Consider the geodesic quadrilateral in $\O$ with
vertices $v'$, $w'$, $w_0$ and $v_0$.
Standard hyperbolic trigonometry (see e.g. \cite[p. 88]{F})
gives that
$$
\aligned
\sinh \frac{d_\O(v',w')}2
& = \sinh \frac{d_\O(v_0,w_0)}2 \, \cosh h(z)
= \cosh h(z) \sqrt{\frac{\cosh d_\O(v_0,w_0) - 1}2}
\\
& = \frac{1}{\sqrt{2}} \, \cosh h(z)
\sqrt{\frac{\tanh d_\O(w,w_0)}{\tanh l_n} -1 \,}  \,
\le \frac{1}{\sqrt{2}} \, \cosh h(z) \sqrt{\frac{1}{\tanh h(z)} -1 }
\\
& = \frac{1}{\sqrt{2}} \, \cosh h(z) \sqrt{\frac{1-\tanh\!^2 h(z)}{\tanh h(z)}}
= \frac{1}{\sqrt{2 \tanh h(z)}}
\le \frac{1}{\sqrt{2 \tanh l_0}} \;.
\endaligned
$$
This fact and Lemma \ref{l:cartesiana} imply
$$
\aligned
d_1(z, v)
& = d_\O(z, v') +  d_\O(v', v)
\le  d_\O(v', w') +  d_\O(z, w') +  d_\O(w', w)
\\
& \le 2 \Arcsinh \frac{1}{\sqrt{2 \tanh l_0}} + d_1(z, w)
\le 2\, d_\O(z, w) + 2 \Arcsinh \frac{1}{\sqrt{2 \tanh l_0}}\;.
\endaligned
$$
\end{proof}

\begin{lema}
\label{l:partirendos}
Let us consider Fermi coordinates $(u,v)$ in $\DD$.
Fix $u_1<u_4$,
$g_1:=\{(u,v):\,u=u_1, \, 0\le v \le x\}$,
$g_4:=\{(u,v):\,u=u_4, \, v\ge 0 \}$,
and $g_2$ the (infinite) geodesic orthogonal to $g_1$ in $(u_1,x)$.
We assume that $g_2$ does not intersects $g_4$.
Consider $(u_4,h)\in g_4$, with $h\ge x$,
and $(u_2,v_2)\in g_2$, with $d\big((u_2,v_2),(u_4,h)\big)=d\big(g_2,(u_4,h)\big)$.
Then
$$
d\big(g_2,(u_4,h)\big)
\le d\big(g_2,(u_3,h)\big) + d\big((u_3,h),(u_4,h)\big)
\le 6 \, d\big(g_2,(u_4,h)\big) \,,
$$
for every $u_2 \le u_3 \le u_4$.
\end{lema}

\begin{proof}
We only need to prove the second inequality.
Fix $u_3 \in [u_2 , u_4]$.

Let us assume that $v_2 \le h$.
Then Lemma \ref{l:cartesiana} implies
$$
\aligned
d\big(g_2,(u_3,h)\big) + d\big((u_3,h),(u_4,h)\big)
& \le d\big((u_2,v_2),(u_2,h)\big)
+ d\big((u_2,h),(u_3,h)\big)
+ d\big((u_3,h),(u_4,h)\big)
\\
& \le d\big((u_2,v_2),(u_2,h)\big)
+ 2 \, d\big((u_2,h),(u_4,h)\big)
\\
& \le 2 \, d_2\big((u_2,v_2),(u_4,h)\big)
\le 6 \, d\big((u_2,v_2),(u_4,h)\big)
= 6 \, d\big(g_2,(u_4,h)\big) \,.
\endaligned
$$
\indent
Let us assume now that $v_2 \ge h$.
Lemma \ref{l:cartesiana} also implies
$$
\aligned
d\big(g_2,(u_3,h)\big) + d\big((u_3,h),(u_4,h)\big)
& \le d\big((u_2,v_2),(u_2,h)\big)
+ d\big((u_2,h),(u_3,h)\big)
+ d\big((u_3,h),(u_4,h)\big)
\\
& \le d\big((u_2,v_2),(u_2,h)\big)
+ 2 \, d\big((u_2,h),(u_4,h)\big)
\\
& \le 2 \, d_1\big((u_2,v_2),(u_4,h)\big)
\le 4 \, d\big((u_2,v_2),(u_4,h)\big)
= 4 \, d\big(g_2,(u_4,h)\big) \,.
\endaligned
$$
\end{proof}

\begin{lema}
\label{l:aex}
Let us define $F$ as
$$
F(a,x) :=
\left\{ \begin{array}{ll}
\displaystyle\frac{1}{\sinh 1}\,\sinh a \cosh x\,, &
\qquad \text{ if } \, 0\le a \le 1\,,
\\
\, & \,
\\
\log \big( \sinh a \cosh x \big) \,, & \qquad \text{ if } \, a\ge 1\,.
\end{array}
\right.
$$
Then
$$
F(a,x) \le a \, e^x \le 2 \sinh a \cosh x \,,
$$
for every $a,x \ge 0$.
\end{lema}

\begin{proof}
The last inequality is a direct consequence of $a \le \sinh a$
and $e^x \le 2 \cosh x$.

If $a\ge 1$, the function $h(x):=a\,e^x-a-x$ satisfies
$h'(x)=a\,e^x-1 \ge a-1 \ge 0$ for every $x \ge 0$.
Hence, $h(x) \ge h(0) = 0$ for every $x \ge 0$,
and we conclude
$$
a \, e^x \ge a+x
= \log \big( e^a e^x \big)
 \ge \log \big( \sinh a \cosh x \big) ,
$$
for $a\ge 1$ and $x \ge 0$.

Since the function $H(a):= \sinh a - a \sinh 1$ is convex in $[0,1]$,
it satisfies
$H(a) \le \max\{H(0),\, H(1)\} = 0$ for every $0\le a \le 1$.
Hence,
$$
a \, e^x \ge  \frac{1}{\sinh 1}\,\sinh a \, e^x
 \ge \frac{1}{\sinh 1}\,\sinh a \cosh x \,,
$$
for $0\le a\le 1$ and $x \ge 0$.
\end{proof}

This result has the following direct corollary.

\begin{corolario}
\label{c:aex}
For a set $E \subset \{(a,x): \, a,x \ge 0\}$, we have
$\Arcsinh \big( \sinh a \cosh x \big) \le c_1$, for every $(a,x) \in E$ and some
constant $c_1$, if and only if
$ a \, e^x \le c_2$, for every $(a,x) \in E$ and some constant $c_2$.

Furthermore, if one of the inequalities holds, the constant in the other
inequality only depends on the first constant.
\end{corolario}

As usual, we denote by $x_+$ the positive part of $x$:
$x_+:=x$ if $x\ge 0$ and $x_+:=0$ if $x< 0$.

\begin{prop}
\label{p:f}
$\,$

\begin{enumerate}
\item[$(1)$] There exists a universal constant $c_1$ such that
$$
f(x,y,t):=\Arccosh \frac{\cosh t + \cosh x \cosh y}{\sinh x \sinh y}
\ge c_1 \big( e^{-x} + e^{-y} + e^{-\frac12 (x+y-t)_+} + (t-x-y)_+ \big) \,,
$$
for every $x,y,t \ge 0$.

\item[$(2)$] For each $l_0>0$, there exists a constant $c_2$,
which only depends on $l_0$, such that
$$
\Arccosh \frac{\cosh t + \cosh x \cosh y}{\sinh x \sinh y}
\le c_2 \big( e^{-x} + e^{-y} + e^{-\frac12 (x+y-t)_+} + (t-x-y)_+ \big) \,,
$$
for every $t \ge 0$ and $x,y \ge l_0$.
\end{enumerate}
\end{prop}

\noindent {\bf Remark.} This result is interesting by itself: if $H$
is a right-angled hexagon in the unit disk for which three pairwise
non-adjacent sides $X$, $Y$, $T$ are given (with respective lengths
$x$, $y$, $t$), then the opposite side of $T$ in $H$ has length
$f(x,y,t)$ (see e.g. \cite[p. 86]{F}, or the proof of Theorem
\ref{t:caractren}).

\begin{proof}
First, we remark that if $x\ge l_0$, then
$e^{-2l_0} e^{2x}\ge 1$ and $e^{2x}-1 \ge (1-e^{-2l_0}) e^{2x}$.
Therefore, if we define $c_3^{-1}:=(1-e^{-2l_0})/2$, we have
$$
e^{2x}-1 \ge 2\, c_3^{-1} e^{2x} \,, \qquad
\sinh x \ge c_3^{-1} e^{x} \,, \qquad
\coth x = 1+ \frac2{e^{2x}-1} \le  1 + c_3\, e^{-2x} \,, \qquad
\text{ for every }
x\ge l_0 \,.
$$
We also have
$$
\coth x = 1+ \frac2{e^{2x}-1} \ge  1 + 2\, e^{-2x} \,, \qquad
\text{ for every }
x\ge 0 \,.
$$

Let us start with the proof of item $(1)$.

If $f \ge 3$, then $f \ge e^{-x} + e^{-y} + e^{-\frac12 (x+y-t)_+}$.
If $f \le 3$, then $1 + \frac23 \,c_4^{-2}f^2 \ge \cosh f$,
for some universal constant $c_4 \le 1$, and
$$
\aligned
1 + \frac23 \,c_4^{-2}f^2
& \ge \cosh f
\ge 2\,e^{t-x-y} + \coth x \coth y
\ge 2\,e^{-(x+y-t)} + \big(1 + 2\, e^{-2x}\big) \big(1 + 2\, e^{-2y}\big) \,,
\\
1 + \frac23 \,c_4^{-2}f^2
& \ge 1 + 2 \big(e^{-2x} + e^{-2y} + e^{-(x+y-t)_+} \big) \,,
\\
c_4^{-1}f
& \ge \sqrt3 \sqrt{e^{-2x} + e^{-2y} + e^{-(x+y-t)_+}}
\ge e^{-x} + e^{-y} + e^{-\frac12 (x+y-t)_+} ,
\\
f & \ge c_4 \big( e^{-x} + e^{-y} + e^{-\frac12 (x+y-t)_+}\big)\,,
\endaligned
$$
where we have used the inequality
$\sqrt3 \sqrt{a + b + c} \ge \sqrt{a} + \sqrt{b} + \sqrt{c}\,$,
for every $a,b,c \ge 0$.
This inequality is $(1)$ if $t \le x+y$.
If $t \ge x+y$, then
$$
\cosh f > \frac{\cosh t}{\sinh x \sinh y} + 1
\ge 2\, e^{t-x-y} + 1
> \frac42 \, e^{t-x-y} + \frac1{4\cdot 2} \, e^{-(t-x-y)}
= \cosh \big( t-x-y + \log 4 \big)
$$
and $f > t-x-y + \log 4 > (t-x-y)_+ + e^{-\frac12 (x+y-t)_+}$.

Consequently we have
$$
f \ge c_1 \big( e^{-x} + e^{-y} + e^{-\frac12 (x+y-t)_+} + (t-x-y)_+ \big) \,,
$$
for every $x,y,t \ge 0$, with $c_1:= c_4/2$, since $c_4 \le 1$.

\spb

Next, let us prove item $(2)$. Fix $l_0>0$. We have seen that $\sinh
x \ge c_3^{-1} e^{x}$ and $\coth x \le  1 + c_3\, e^{-2x}$, for
every $x\ge l_0$.

Let us assume $t\ge x+y$.
If $x,y \ge l_0$, then
$$
\frac12 \,e^{f}
\le \cosh f
= \frac{\cosh t + \cosh x \cosh y}{\sinh x \sinh y}
\le c_3^2 \, e^{t-x-y} + \cotanh\!^2 l_0 \,.
$$
Consequently,
$$
e^{f}
\le 2\, c_3^2 \, e^{t-x-y} + 2 \cotanh\!^2 l_0
\le e^{t-x-y+c_5} ,
$$
with $c_5:= \log \big(2\, c_3^2 + 2 \cotanh\!^2 l_0\big)$,
since $t- x-y \ge 0$.
Hence,
$f \le t-x-y +c_5 = (t-x-y)_+ + c_5 e^{-\frac12 (x+y-t)_+}$,
for every $t \ge 0$ and $x,y \ge l_0$ with $t\ge x+y$.

Let us assume $t\le x+y$.
If $x,y \ge l_0$, then
$$
\aligned
1 + \frac12 \,f^2
& \le \cosh f
\le c_3^2 \, e^{t-x-y} + \cotanh x \cotanh y
\le c_3^2 \, e^{t-x-y} + \big(1 + c_3\, e^{-2x}\big) \big(1 + c_3\, e^{-2y}\big) \,,
\\
\frac12 \,f^2
& \le c_3^2 \, e^{t-x-y} + c_3\, e^{-2x} + c_3\, e^{-2y} + c_3^2\, e^{-2x-2y} ,
\\
\frac12 \,f^2
& \le c_3^2 \, e^{t-x-y} + c_3\, e^{-2x} + c_3\, e^{-2y} + \frac12\,
c_3^2\, \big( e^{-2x} + e^{-2y}\big) \,,
\\
f^2
& \le 2 \, c_3^2 \, e^{-(x-y-t)} + \big(2\,c_3 + c_3^2\big) e^{-2x}
+ \big(2\,c_3 + c_3^2\big) e^{-2y} ,
\\
f^2
& \le c_6^2 \big( e^{-2x} + e^{-2y} + e^{-(x+y-t)_+} \big) \,,
\\
f
& \le c_6 \big( e^{-x} + e^{-y} + e^{-(x+y-t)_+} + (t-x-y)_+ \big) \,,
\endaligned
$$
where $c_6^2:= \max \big\{ 2\,c_3^2 ,\, 2\, c_3 + c_3^2 \big\}$,
for every $t \ge 0$ and $x,y \ge l_0$ with $t\le x+y$.
Then we have $(2)$ with $c_2:=\max\{1, \, c_5, \, c_6\}$.
\end{proof}

\begin{prop}
\label{p:F}
For each $l_0>0$, we have
$$
F(x,y,t,h):=\Arcsinh \frac{\cosh x \cosh (y-h) + \cosh t \cosh h}{\sinh y}
\asymp  e^{-h+x} + e^{-(y-h-t)_+} + (t+h-y)_+ \,,
$$
for every $x,y,t,h \ge 0$, verifying $y \ge h \ge x$ and $y \ge l_0$.
Furthermore, the constants in the inequalities only depend on $l_0$.
\end{prop}

\noindent {\bf Remark.} This result is interesting by itself: if $H$
is a right-angled hexagon in the unit disk for which three pairwise
non-adjacent sides $X$, $Y$, $T$ are given (with respective lengths
$x$, $y$, $t$), $P$ is the nearest point to $X$ in $Y$, and $P_h$ is
the point in $Y$ with $d(P_h,P)=h$, then $F(x,y,t,h)$ is the
distance between $P_h$ and the opposite side of $Y$ in $H$ (see the
proof of Theorem \ref{t:caractren}).

\begin{proof}
We have seen that if $y\ge l_0$, and $c_3^{-1}:=(1-e^{-2l_0})/2$, we have
$c_3^{-1} e^{y}\le \sinh y \le e^{y}/2$.
We also have
$e^{z}/2\le \cosh z \le e^{z}$, for every $z\ge 0$.

Then $\sinh F \asymp e^{-h+x} + e^{-y+h+t}$,
since $y \ge l_0$ and $y \ge h$,
and the constants in the inequalities only depend on $l_0$.

If $h+t \le y$, then $e^{-h+x} + e^{-y+h+t} \le 2$, and
$$
F \asymp \sinh F \asymp e^{-h+x} + e^{-(y-h-t)}
=  e^{-h+x} + e^{-(y-h-t)_+} + (t+h-y)_+ \,.
$$
\indent
If $h+t \ge y$, then $e^{-h+x} + e^{-y+h+t} \ge 1$, and
$$
e^F \asymp \sinh F \asymp e^{-h+x} + e^{-y+h+t}
\asymp e^{t+h-y} = e^{-1} e^{1+ (t+h-y)_+} \,.
$$
Since
$$
F \ge \Arcsinh \frac{\big(e^x e^{y-h} + e^t e^h\big)/4}{e^{y}/2}
\ge \Arcsinh \frac12 \big( e^{-h+x} + e^{-y+h+t}\big)
\ge \Arcsinh \frac12>0 \,,
$$
and $1+ (t+h-y)_+ \ge 1>0$
for every $x,y,t,h \ge 0$,
and $e^F \asymp e^{1+ (t+h-y)_+}$
for every $x,y,t,h \ge 0$, verifying $h+t \ge y \ge h \ge x$ and $y \ge l_0$,
we obtain that
$F \asymp 1+ (t+h-y)_+$.
Since $1\le e^{-h+x} + 1 = e^{-h+x} + e^{-(y-h-t)_+} \le 2$,
we also conclude that
$F \asymp e^{-h+x} + e^{-(y-h-t)_+} + (t+h-y)_+$, if $h+t \ge y$.
\end{proof}

The following corollary can be directly deduced from this result.

\begin{corolario}
\label{c:xyth}
For each $l_0>0$, let us consider
a set $E \subset \{(x,y,t,h): \, x,y,t,h \ge 0, \,
y \ge h \ge x, \, y \ge l_0 \}$.
We have
$F(x,y,t,h) \le c_1$, for every $(x,y,t,h) \in E$ and some constant $c_1$,
if and only if
$(t+h-y)_+ \le c_2$, for every $(x,y,t,h) \in E$ and some constant $c_2$.

Furthermore, if one of the inequalities holds, the constant in the other
inequality only depends on the first constant and $l_0$.
\end{corolario}

Obviously, we can replace condition
$(t+h-y)_+ \le c_2$ by $t+h-y \le c_2$.
We prefer the first one since $F$ will be a distance and
$(t+h-y)_+\ge 0$.




\


\begin{thebibliography}{99}









\bibitem{A} Aikawa, H.,
Positive harmonic functions of finite order in a Denjoy type domain,
{\it Proc. Amer. Math. Soc.} {\bf 131} (2003), 3873-3881.


\bibitem{APR} Alvarez, V., Pestana, D., Rodr{\'\i}guez, J. M.,
Isoperimetric inequalities in Riemann surfaces of infinite type,
{\it Rev. Mat. Iberoamericana} {\bf 15} (1999), 353-427.


\bibitem{APRT} Alvarez, V., Portilla, A., Rodr{\'\i}guez, J. M., Tour{\'\i}s, E.,
Gromov hyperbolicity of Denjoy domains,
{\it Geom. Dedicata} {\bf 121} (2006), 221-245.



\bibitem{ARY} Alvarez, V., Rodr\'\i guez, J.M., Yakubovich, D.V.,
Subadditivity of p-harmonic ``measure" on graphs,
{\it Michigan Math. J.} {\bf 49} (2001), 47-64.



\bibitem{BB} Balogh, Z. M., Buckley, S. M.,
Geometric characterizations of Gromov hyperbolicity,
{\it Invent. Math.} {\bf 153} (2003), 261-301.

\bibitem{B1} Basmajian, A.,
Constructing pair of pants,
{\it Ann. Acad. Sci. Fenn. Series AI} {\bf 15} (1990), 65-74.

\bibitem{B2} Basmajian, A.,
Hyperbolic structures for surfaces of infinite type,
{\it Trans. Amer. Math. Soc.} {\bf 336} (1993), 421-444.



\bibitem{Be} Benoist, Y.,
Convexes hyperboliques et fonctions quasisymtriques,
{\it Publ. Math. Inst. Hautes \'Etudes Sci.} {\bf 97} (2003), 181-237.

\bibitem{BHK} Bonk, M., Heinonen, J., Koskela, P.,
Uniformizing Gromov hyperbolic spaces. Ast\'erisque No. 270 (2001).





\bibitem{CFPR} Cant\'on, A., Fern\'andez, J. L., Pestana, D., Rodr\'{\i}guez, J. M.,
On harmonic functions on trees,
{\it Potential Analysis} {\bf 15} (2001), 199-244.

\bibitem{C} Chavel, I., {\it Eigenvalues in Riemannian Geometry.}
Academic Press, New York, 1984.






\bibitem{F} Fenchel, W., {\it Elementary Geometry in Hyperbolic Space.}
Walter de Gruyter, Berlin-New York, 1989.


\bibitem{FR2} Fern\'andez, J. L., Rodr\'{\i}guez, J. M.,
Area growth and Green's function of Riemann surfaces,
{\it Arkiv f\"or matematik} {\bf 30} (1992), 83-92.


\bibitem{GJ} Garnett, J., Jones, P.,
The Corona theorem for Denjoy domains,
{\it Acta Math.} {\bf 155} (1985), 27-40.


\bibitem{GH} Ghys, E., de la Harpe, P.,
Sur les Groupes Hyperboliques d'apr\`es Mikhael Gromov.
Progress in Mathematics, Volume 83. Birkh\"auser. 1990.


\bibitem{G} Gonz\'alez, M. J.,
An estimate on the distortion of the logarithmic capacity,
{\it Proc. Amer. Math. Soc.} {\bf 126} (1998), 1429-1431.



\bibitem{G1} Gromov, M.,
Hyperbolic groups, in ``Essays in group theory".
Edited by S. M. Gersten, M. S. R. I. Publ. {\bf 8}. Springer, 1987, 75-263.


\bibitem{G2} Gromov, M.
(with appendices by M. Katz, P. Pansu and S. Semmes),
Metric Structures for Riemannian and Non-Riemannnian Spaces.
Progress in Mathematics, vol. 152. Birkh\"auser, 1999.


\bibitem{H} Haas, A.,
Dirichlet points, Garnett points and infinite ends of hyperbolic surfaces I,
{\it Ann. Acad. Sci. Fenn. Series AI} {\bf 21} (1996), 3-29.

\bibitem{Ha} H\"ast\"o, P. A.,
Gromov hyperbolicity of the $j_G$ and $\widetilde{j}_G$ metrics,
{\it Proc. Amer. Math. Soc.} In press.







\bibitem{HS} Holopainen, I., Soardi, P. M., $p$-harmonic functions on
graphs and manifolds,
{\it Manuscripta Math.} {\bf 94} (1997), 95-110.


\bibitem{K1} Kanai, M., Rough isometries and combinatorial approximations
of geometries of non-compact Riemannian manifolds, {\it J. Math. Soc.
Japan} {\bf 37} (1985), 391-413.

\bibitem{K2} Kanai, M., Rough isometries and the parabolicity of
Riemannian manifolds, {\it J. Math. Soc. Japan} {\bf 38} (1986),
227-238.

\bibitem{K3} Kanai, M., Analytic inequalities and rough isometries
between non-compact Riemannian manifolds.
Curvature and Topology of Riemannian manifolds (Katata, 1985).
Lecture Notes in Math. {\bf 1201}. Springer (1986), 122-137.

\bibitem{KN} Karlsson, A., Noskov, G. A.,
The Hilbert metric and Gromov hyperbolicity,
{\it Enseign. Math.} {\bf 48} (2002), 73-89.







\bibitem{PRT1} Portilla, A., Rodr\'{\i}guez, J. M., Tour{\'\i}s, E.,
Gromov hyperbolicity through decomposition of metric spaces II,
{\it J. Geom. Anal.}  {\bf 14} (2004), 123-149.

\bibitem{PRT2} Portilla, A., Rodr\'{\i}guez, J. M., Tour{\'\i}s, E.,
The topology of balls and Gromov hyperbolicity of Riemann surfaces,
{\it Diff. Geom. Appl.} {\bf 21} (2004), 317-335.

\bibitem{PRT3} Portilla, A., Rodr\'{\i}guez, J. M., Tour{\'\i}s, E.,
The role of funnels and punctures in the Gromov hyperbolicity of Riemann surfaces.
{\it Proc. Edinburgh Math. Soc.} {\bf 49} (2006), 399-425.

\bibitem{PT} Portilla, A., Tour{\'\i}s, E.,
A characterization of Gromov hyperbolicity of surfaces
with variable negative curvature,
{\it Publicacions Matem\`atiques}, to appear.








\bibitem{R1} Rodr\'{\i}guez, J. M.,
Isoperimetric inequalities and Dirichlet functions of Riemann surfaces,
{\it Publicacions Matem\`atiques} {\bf 38} (1994), 243-253.

\bibitem{R2} Rodr\'{\i}guez, J. M.,
Two remarks on Riemann surfaces,
{\it Publicacions Matem\`atiques} {\bf 38} (1994), 463-477.

\bibitem{RT1} Rodr\'{\i}guez, J. M., Tour{\'\i}s, E.,
Gromov hyperbolicity through decomposition of metric spaces,
{\it Acta Math. Hung.} {\bf 103} (2004), 53-84.

\bibitem{RT2} Rodr\'{\i}guez, J. M., Tour{\'\i}s, E.,
A new characterization of Gromov hyperbolicity for Riemann surfaces,
{\it Publicacions Matem\`atiques} {\bf 50} (2006), 249-278.

\bibitem{RT3} Rodr\'{\i}guez, J. M., Tour{\'\i}s, E.,
Gromov hyperbolicity of Riemann surfaces,
{\it Acta Math. Sinica} {\bf 23} (2007), 209-228.







\bibitem{So} Soardi, P. M., Rough isometries and Dirichlet finite harmonic
functions on graphs, {\it Proc. Amer. Math. Soc.} {\bf 119} (1993),
1239-1248.







\bibitem{V} V\"ais\"al\"a, J., Hyperbolic and uniform domains in Banach spaces.
Preprint.  \newline
[Available at \email{http://www.helsinki.fi/\~\,  jvaisala/preprints.html}]

\end{thebibliography}
\end{document}